
\ifx\CWIloaded\undefined \let\CWIloaded\relax \else   \fi

\ifx\topglue\undefined \def\topglue{\nointerlineskip\vglue-\topskip\vskip}\fi

\ifx\formloaded\undefined 
  \vsize=24.5 cm \advance \vsize -24pt 
  \hsize=16 cm
  \let\formloaded\relax 
\fi


\catcode`_=11 

\def\let_next{\let\next=}
{\obeylines\gdef\default_address%
 {Universit\'e de Poitiers, D\'epartement de Math\'ematiques,
  SP2MI, T\'el\'eport~2, BP~179, \ 86960 Futuroscope Cedex, France}}

\newtoks\title \newtoks\MSCcodes \newtoks\CRcodes \newtoks\keywords
\newtoks\titlenote
\newif\if_started

 3

\font\seventeenssbf=cmssbx10 scaled \magstep 3

\font\twelvess=cmss12

\font\tenss=cmss10
\font\tenssit=cmssi10
\font\tenssbf=cmssbx10
\def\sans
{\def\rm{\fam0\tenss}\def\it{\fam\itfam\tenssit}\def\bf{\fam\bffam\tenssbf}\rm
}


\font\ninerm=cmr9 \font\nineit=cmti9 \font\ninebf=cmbx9 \font\ninesl=cmsl9
\font\ninett=cmtt9 \font\ninei=cmmi9 \font\ninesy=cmsy9
\font\niness=cmss9 \font\ninessit=cmssi9 \font\ninessbf=cmssbx10 at 9 pt
\skewchar\ninei='177 \skewchar\ninesy='60

\def\ninept
{\def\rm{\fam0\ninerm}\def\it{\fam\itfam\nineit}\def\bf{\fam\bffam\ninebf}%
 \def\sl{\fam\slfam\ninesl}\def\tt{\fam\ttfam\ninett}\rm
 \def\sans{\def\rm{\fam0\niness}\def\it{\fam\itfam\ninessit}%
           \def\sl{\fam\slfam\ninessit}\def\bf{\fam\bffam\ninessbf}\rm}%
 \textfont0\ninerm \textfont1\ninei \textfont2=\ninesy
 \textfont\itfam\nineit \textfont\slfam\ninesl
 \textfont\bffam\ninebf \textfont\ttfam\ninett
 \normalbaselineskip=11pt \normalbaselines
 \setbox\strutbox=\hbox{\vrule height 8 pt depth 3 pt width 0 pt}%
}


\long\def\centerpar#1{\begingroup
  \advance\leftskip by 0pt plus 1 fil minus 3pt \rightskip=\leftskip
	\parskip=0pt \parindent=0pt \parfillskip=0pt \linepenalty=100
  \noindent\ignorespaces#1\par\endgroup}

\def\good_break{\unskip\penalty-10\ \ignorespaces}

\catcode`@=11
\def\vfootnote#1%
{\insert\footins\bgroup \ninept
 \interlinepenalty=\interfootnotelinepenalty
 \splittopskip=\ht\strutbox 
 \splitmaxdepth=\dp\strutbox 
 \floatingpenalty=20000 
 \leftskip=1pc \rightskip=1pc \spaceskip=0pt \xspaceskip=0pt
 \parindent=0pt \parskip=3pt minus 1pt \parfillskip=0pt plus 1 fil
 \textindent{#1}\footstrut\futurelet\next\f@@t
}
\catcode`@=12

\let\_authors\empty 

%

\def\Author#1{\def\this_author{#1}%
  \ifx\this_author\_eq \let_next\Author 
  \else \let_next\let_next \afterassignment\_author 
  \fi \next}
\def\_author
 {\let\this_addr=\default_address
  \ifx\next\address \let_next\get_address
  \else \ifx\next\email \let_next\get_email
  \else \ifx\next\URL \let_next\get_URL
  \else \add_author
  \fi\fi\fi \next
}

\def\address#1{\errmessage
 {\string\address\space should only follow \string\Author}}
\def\email#1{\errmessage
 {\string\email\space should only follow \string\Author}}
\def\URL#1{\errmessage
 {\string\URL\space should only follow \string\Author}}

\def\get_address{\begingroup\obeylines\get_addr}
\def\get_addr#1%
{\endgroup\def\this_addr{#1}%
  \ifx\this_addr\_cr \let_next\get_address 
  \else\ifx\this_addr\_eq \let_next\get_address 
       \else \let_next\let_next \afterassignment\get_ad
  \fi\fi \next
}
\def\get_ad
{\ifx\next\email \let_next\get_email
 \else \ifx\next\URL \let_next\get_URL
 \else \add_author
 \fi\fi \next
}

\def\get_email#1%
{\def\this_email{#1}%
  \ifx\this_email\_eq \let_next\get_email
  \else \edef\this_addr{\this_addr\endgraf{\tt #1}}%
    \let_next\let_next \afterassignment\after_em
  \fi \next
}
\def\after_em{\ifx\next\URL \let_next\get_URL \else \add_author \fi \next}
\def\get_URL#1%
{\def\this_URL{#1}%
 \ifx\this_URL\_eq \let_next\get_URL
 \else 
  {\def~{\char "7E }%
   \xdef\this_addr{\this_addr\endgraf{\noexpand\tt URL: http://#1}}%
  }\let_next\add_author
 \fi \next
}

\def\_eq{=}{\obeylines\gdef\_cr{^^M}}

\def\add_author 
{\toks0=\expandafter{\_authors}%
 \xdef\_authors
 {\the\toks0 \noexpand\\{\this_author\noexpand\_addr{\this_addr}}}%
}


\def\start_art
{\vfil\eject \ifodd\pageno\else \shipout\vbox to \vsize{}\advancepageno\fi

 \topglue 10mm plus 3 mm minus 2mm
 \centerpar{\let\\=\good_break\seventeenssbf \baselineskip=22pt \the\title}%
 \ifx \_authors\empty \message{Warning: no authors specified.}\else
   \bigskip\centerpar
    {\twelvess
     \def\par{\ifhmode\endgraf\else\medskip\fi}\obeylines
     \def\_addr##1{\smallskip{\ninessit\let\\=\par\noindent##1\endgraf}}%
     \def\\##1{\bigskip##1}\_authors
    }\bigskip
 \fi
 \_startedtrue
}


\outer\def\abstract
{\if_started \errmessage{Abstract must be placed at start of article.}%
 \else \start_art \bigskip \begingroup \ninept \sans \narrower
   \leftline{\hskip\leftskip ABSTRACT}\nobreak\smallskip
	 \everypar{\everypar{}\setbox0=\lastbox}%
 \fi
}

\outer\def\maintext
{\if_started\medskip\endgroup\else\start_art\fi
 \if &\the\MSCcodes\the\CRcodes\the\keywords\the\titlenote&\else
  {\narrower\ninept\sans
   \if &\the\MSCcodes&\else
     \hangindent=4em \noindent
     {\it 1991 Mathematics Subject Classification:\/} \the\MSCcodes.
     \par\fi
   \if &\the\CRcodes&\else
     \hangindent=4em \noindent
     {\it 1991 Computing Reviews Subject Classification:\/} \the\CRcodes.
     \par\fi
   \if &\the\keywords&\else
     \hangindent=2em \noindent{\it Keywords and Phrases:\/} \the\keywords.
     \par\fi
   \if &\the\titlenote&\else
     {\def\\{\unskip\hfil\break\ignorespaces}%
      \hangindent=2em \noindent{\it Note:\/} \the\titlenote.
      \par}\fi
  }\bigskip
 \fi
}

\def\sticksection 
{\vbox\bgroup\parskip=0pt \everypar{\egroup\noindent}}

\def\newreport
{\vfil \eject \secno=0\subsecno=0\proclno=0 \eqnumber=0 \pageno=1 \mark{}
 \title={}\MSCcodes={}\CRcodes={}\keywords={}\titlenote={}
 \let\_authors\empty
 \_startedfalse
}

\catcode`_=8 

\input marxmax

\refindent=3pc
\normalbaselineskip=13pt\normallineskiplimit=.1pt\normallineskip.1pt
\normalbaselines
\setbox\strutbox=\hbox{\vrule height8.9pt depth4pt width0pt}


\def\RS.{Robinson-Schensted}
\def\Sch.{Sch\"utzen\-berger}

\def\al.{algorithm}
\def\ta.{tableau}
\def\irr.{irreducible}
\def\jdt.{jeu de taquin}
\def\Jdt.{Jeu de taquin}

\opdef\sh  \opdef\ch
\opdef\im
\altopdef\Im{im}
\altopdef\Ker{ker}
\def\tr{^{\rm t}}

\mathorddef\D=D2 
\mathorddef\F=F2 
\mathorddef\E=E2 
\mathorddef\G=G2 
\mathorddef\U=U2 
\mathorddef\Part=P2 
\mathorddef\Tab=T2 
\mathorddef\B=B\bffam
\mathorddef\Sym=S\bffam
\def\GL{{\bf GL}}
\def\S#1{\Sym_{#1}}
\def\Sn{\Sym_n}
\def\w{\tilde n}

\def\[#1]{_{/#1}}
\def\Fu{\F_\eta}
\def\FuH{\F_{\eta|_H}}
\def\Ful{\F_{\eta\[l]}}
\def\au{\alpha_\eta}
\def\ou{\omega_\eta}
\def\Tx{\Tab_{\\/\mu}^{\,\slid\nu}}
\def\Gr#1{\G_{#1}^\eta}
\def\Gx{\G^\eta_{\mu,\nu}}
\def\CT#1{{\bf T}_{#1}}

\mathorddef\nil=\odot
\mathdef\adj=Rel \parallel
\def\nadj{\mathrel{\rlap{\kern-.75pt\it/}}\adj}
\def\ar(#1,#2,#3,#4){\smash{\({#1\atop#3}\,{#2\atop#4}\)}}
\def\Ar(#1,#2,#3,#4){\pmatrix{#1&#2\cr#3&#4\cr}}
\mathdef\slid=Rel \triangleright

{\catcode`\"=\active \gdef\letqq{\let"=}}
{\catcode`\@=11
 \def\next{^\bgroup\d@wn}\global\letqq\next
 \gdef\d@wn{\downarrow\futurelet\next\pr@m@s}
 \gdef\pr@m@s{\ifx'\next\letnext\pr@@@s \else\ifx^\next\letnext\pr@@@t
              \else\ifx"\next\letnext\d@@n\else\ifx!\next\letnext\d@@@n
              \else\letnext\egroup\fi\fi\fi\fi
              \next}
 \gdef\d@@n"{\d@wn}\gdef\d@@@n!{\d@wn}
}
\mathcode`\"="8000 

\def\hi#1{\lceil#1\rceil}
\def\lo#1{\lfloor#1\rfloor}
\input \ta.x

\subsecno=-1 

\readrefs

\Author={Marc A. A. van Leeuwen}
\email{maavl@mathlabo.univ-poitiers.fr}
\URL{wwwmathlabo.univ-poitiers.fr/~maavl/}

\title={Flag Varieties and \\ Interpretations of Young Tableau Algorithms}

\MSCcodes{05E15, 20G15}


\titlenote
{Work partly supported by the {\it Schweizerische Nationalfonds} and by
 {\it l'Univerist\'e du Qu\'ebec \`a Montr\'eal}}

\abstract

The conjugacy classes of nilpotent $n\times n$ matrices can be parametrised by
partitions~$\\$ of~$n$, and for a nilpotent~$\eta$ in the class parametrised
by~$\\$, the variety~$\Fu$ of $\eta$-stable flags has its irreducible
components parametrised by the standard Young tableaux of shape~$\\$. We
indicate how several algorithmic constructions defined for Young tableaux have
significance in this context, thus extending Steinberg's result that the
relative position of flags generically chosen in the irreducible components
of~$\Fu$ parametrised by tableaux $P$~and~$Q$, is the permutation associated
to $(P,Q)$ under the Robinson-Schensted correspondence. Other constructions
for which we give interpretations are \Sch.'s involution of the set of Young
tableaux, jeu de taquin (leading also to an interpretation of
Littlewood-Richardson coefficients), and the transpose Robinson-Schensted
correspondence (defined using column insertion). In each case we use a doubly
indexed family of partitions, defined in terms of the flag (or pair of flags)
determined by a point chosen in the variety under consideration, and we show
that for generic choices, the family satisfies combinatorial relations that
make it correspond to an instance of the algorithmic operation being
interpreted, as described in~\ref{festschrift}.

\maintext

\secno=-1
\vbox{
\newsection Introduction.

The Schensted algorithm, which defines a bijective correspondence between
permutations and pairs of (standard) Young tableaux, and the \Sch. (or
evacuation) algorithm, which defines a shape preserving involution of the set
of Young tableaux, can both be described using doubly indexed families of
partitions that satisfy certain local rules, as described
in~\ref{festschrift}. In this paper we show how both correspondences occur in
relation to questions concerning varieties of flags stabilised by a fixed
nilpotent transformation~$\eta$. The mentioned doubly indexed families of
partitions arise very naturally in this context, and they provide detailed
information concerning the internal steps of the algorithms, rather than just
about the correspondences defined by them. As a consequence, the study of the
\Sch. algorithm also leads to an interpretation of \foreign{jeu de taquin},
and of the Littlewood-Richardson coefficients. The connections between
geometry and combinatorics presented here include and extend results due to
Steinberg (\ref{SteinbergDUV},~\ref{SteinbergORS}) and
Hesselink~(\ref{Hesselink}).

}

The basic fact underlying these interpretations is that the irreducible
components of the variety~$\Fu$ of $\eta$-stable complete flags is
parametrised in a natural way by the set of standard Young tableaux of shape
equal to the Jordan type~$J(\eta)$ of~$\eta$. In fact there are two dual
parametrisations, and we show that the transition between them is given by the
\Sch. involution. Taking a projection on varieties of incomplete flags by
forgetting the parts of flags below a certain dimension, we obtain from this
an interpretation of
jeu de taquin (operating on skew standard tableaux), and a bijection between
Littlewood-Richardson tableaux and irreducible components of a variety of
$\eta$-stable subspaces of fixed type and cotype.

The other interpretations involve the Robinson-Schensted algorithm and
relative positions of flags. We give a derivation of Steinberg's result that
the relative position of two generically chosen flags in given irreducible
components of~$\Fu$ is the permutation related by the Robinson-Schensted
correspondence to the pair of the tableaux parametrising those components. For
this result a specific choice of parametrisation is required, but there are
variations of the statement for other choices; in particular when dual
parametrisations are used for the two components, one obtains the transpose
Robinson-Schensted correspondence, defined using column insertion. Together,
these interpretations give a geometric meaning to many key properties of the
combinatorial correspondences: involutivity of the \Sch. correspondence,
confluence of jeu de taquin, the fact that the number of skew tableaux of a
fixed shape and given rectification~$P$ depends only on the shape of~$P$,
symmetry of the Robinson-Schensted correspondence, and the various relations
between this correspondence, its transpose, and the \Sch. correspondence.

This paper is organised as follows. In \Sec1 we review the essential facts in
the linear algebra of a vector space equipped with a nilpotent
transformation~$\eta$, followed in \Sec2 by the definition of the
variety~$\Fu$, and the parametrisations of its irreducible components by Young
tableaux. In \Sec3 and \Sec4 we give the respective interpretations of the
\Sch. algorithm, and of jeu de taquin and Littlewood-Richardson tableaux. In
\Sec5 we discuss relative positions of flags, which are used in \Sec6 to give
a geometric interpretation of the Robinson-Schensted correspondence, and in
\Sec7 to give a similar but independent interpretation of the transposed
Robinson-Schensted correspondence.

In our combinatorial notation, as well as in our approach of the algorithms
considered, we shall closely follow~\ref{festschrift}; we collect here the
essential definitions used. A \newterm{partition}~$\\$ is an infinite weakly
decreasing sequence $\\_0\geq\\_1\geq\cdots$ of natural numbers (called the
\newterm{parts} of~$\\$) with finite sum, denoted by~$|\\|$. To each
partition~$\\$ is associated its \newterm{Young diagram}
$Y(\\)\ssubset\N\times\N$, defined by $(i,j)\in Y(\\)\iff j<\\_i$ (so that
$\Card{Y(\\)}=|\\|$); the \newterm{transposed} partition~$\\\tr$ of~$\\$ is
the one whose Young diagram is obtained by reflection of~$Y(\\)$ in the main
diagonal. The elements of~$Y(\\)$ are called its \newterm{squares}, and are
depicted correspondingly (so that they may be filled with values); Young
diagrams are displayed with the first index~$i$ increasing downwards, the
second increasing towards the right, like matrices. The set of all partitions
is denoted by~$\Part$, and is partially ordered by inclusion of Young
diagrams, written~`$\subset$'; the elements of
$\Part_n=\setof\\\in\Part:|\\|=n\endset$, are called partitions of~$n$. For
the predecessor relation in~$\Part$, that was denoted by $\mu\in\\^-$ in
\ref{festschrift}, we shall write instead $\mu\prec\\$, with $\mu\preceq\\$
meaning that $\mu\prec\\$ or $\mu=\\$. When $\mu\prec\\$ and $Y(\\)\setminus
Y(\mu)=\{x\}$ we call the square~$x$ a corner of~$\\$ and a cocorner of~$\mu$,
and write $\\=\mu+x$, \ $\mu=\\-x$ and $\\-\mu=x$. We write $x\adj y$ to
indicate that squares $x$~and~$y$ are adjacent.

A \newterm{Young tableau} is an injective map $T\:Y(\\)\to\N$, for
some~$\\\in\Part$ called the \newterm{shape} $\sh T$ of~$T$, such that when
each number~$T(i,j)$ is written as entry into square~$(i,j)$, all rows and
columns are increasing. Each such~$T$ determines a saturated decreasing chain
$\ch T$ from $\\$~to~$(0)$ in~$\Part$, by recording the successive shapes as
the squares are removed from~$Y(\\)$ in order of decreasing entries. When $\ch
T=\ch T'$ we write $T\sim T'$, which is an equivalence relation on tableaux;
for tableaux of fixed shape~$\\$, a set of representatives of the equivalence
classes is formed by the set $\Tab_\\$ of \newterm{normalised} tableaux, whose
entries are all~$<|\\|$. For any non-empty Young tableau~$T$, the square
containing the highest entry of~$T$ is denoted by~$\hi T$, and the tableau
obtained by removing that square (and its entry) from~$T$ by~$T^-$. By
applying the deflation procedure used in the definition of the \Sch. algorithm
to~$T$, a tableau~$T"$ is obtained, in which the smallest entry has
disappeared. The result of applying the full \Sch. algorithm to~$T$ is denoted
by~$S(T)$, and the pair of tableaux obtained by applying the Schensted
algorithm to a permutation~$\sigma$ is denoted by~$RS(\sigma)$ (see
\ref{festschrift} for definitions).

When numbering or indexing with natural numbers, these definitions all start
using~$0$, rather than~$1$; this holds in particular for the parts of a
partition, and the rows and columns of Young diagrams. This leads to simpler
expressions, but note that while row~$i$ of $Y(\\)$ has length~$\\_i$, it has
no square in column~$\\_i$. We even have gone a bit further than
in~\ref{festschrift}, by defining entries of normalised tableaux to start
at~$0$. Also, the group $\Sn$ is taken to consist of permutations of
$\{0,\ldots,n-1\}$, represented by sequences of the form
$\sigma_0,\sigma_1,\ldots,\sigma_{n-1}$; in particular, the order reversing
permutation~$\w\in\Sn$ satisfies $\w_i=n-1-i$.

The work presented here was motivated by that of the author's
thesis~\ref{Marc}, which deals with the significantly more complicated case of
other classical groups (in characteristic~$\neq2$) instead of~$\GL_n$. There a
combinatorial algorithm analogous to the Robinson-Schensted algorithm is
deduced, that perfoms the corresponding computation of generic relative
positions of flags; however there are complications that cause the
descriptions and proofs to be much more technical than those in the current
paper. Our aim here was to separate the main techniques and arguments
of~\ref{Marc} from a number of distracting technicalities, by applying them in
the simpler situation of~$\GL_n$; even so the reasoning is sometimes detailed
and subtle. At the same time we believe that in describing that situation as
transparently as possible, we have been able to provide some new insight.
 \par
\newsection Nilpotent transformations.

Let $V$ be a vector space of finite dimension~$n$ over an infinite field~$k$.
In this section we recall some basic facts concerning $V$, equipped with a
fixed nilpotent endomorphism~$\eta$.
Most of these facts are also discussed, from a somewhat more general
and elevated perspective, in \ref{Macdonald, Chapter~II}.

A~subspace $V'$ of~$V$ is called $\eta$-stable if $\eta(V')\subset V'$. There
exists a decomposition of~$V$ as a direct sum of non-zero $\eta$-stable
subspaces that cannot be so decomposed further; any summand of dimension~$d$
admits a basis $\li(x0..d-1)$ such that $\eta(x_0)=0$ and $\eta(x_i)=x_{i-1}$
for $0<i<d$. This is called a decomposition into \newterm{Jordan blocks}; it
is generally not unique, but the multiset of the dimensions of the blocks
depends only on~$\eta$. Arranged into weakly decreasing order these dimensions
form a partition of~$n$, called the \newterm{Jordan type} $J(\eta)$ of~$\eta$.
Throughout this paper we write $\\$ for~$J(\eta)$ and $u$ for the unipotent
transformation $\eta+\Id\in\GL(V)$ corresponding to~$\eta$; when $\eta$ is
variable, $\\$~and~$u$ are assumed to vary correspondingly.

One can characterise $\\$ in terms of the powers $\eta^j$ of $\eta$ (among
which we include $\eta^0=\1$), without referring to any particular
decomposition into Jordan blocks, as follows.

\proclaim Proposition. \Jordanprop
For all $c\geq0$ one has $\dim(\Ker\eta^c)=\sum_{0\leq j<c}\\\tr_j$, which is
the number of squares in the first~$c$ columns of the Young diagram $Y(\\)$.
Similarly, $\dim(\Im\eta^c)=\sum_{j\geq c}\\\tr_j$, which is the number of
squares in the remaining columns of~$Y(\\)$.

\proof
This is easily verified for individual Jordan blocks, from which the general
case follows.
\QED

For any $\eta$-stable subspace~$V'$ of~$V$, $\eta$ induces nilpotent
endomorphisms of the spaces $V'$ and $V/V'$, which will be denoted
respectively by $\eta|_{V'}$~and~$\eta\[V']$. Since a nilpotent endomorphism
of a $1$-dimensional space is necessarily~$0$, a subspace~$l$ of dimension~$1$
(a line) is $\eta$-stable if and only if $l\subset\ker\eta$, and similarly a
subspace~$H$ of codimension~$1$ (a hyperplane) is $\eta$-stable if and only if
$H\supset\im\eta$.

\proclaim Proposition. \Jsubsetprop
Let $V'$ be an $\eta$-stable subspace of~$V$.
Then $J(\eta|_{V'})\subset\\$, and $J(\eta\[V'])\subset\\$.

\proof
The image of the subspace $\ker\eta^j$ under the projection $V\to V/V'$ is
$\ker\eta^j/\ker(\eta|_{V'})^j$, for all $j\geq0$; since
$\ker\eta^j\subset\ker\eta^{j+1}$, it follows that
$\dim(\ker\eta^j)-\dim(\ker(\eta|_{V'})^j)$ increases weakly as $j$ increases:
$$\dim(\ker\eta^j)-\dim(\ker(\eta|_{V'})^j)\leq
 \dim(\ker\eta^{j+1})-\dim(\ker(\eta|_{V'})^{j+1})
.\nn
$$
The length $\\\tr_j$ of column~$j$ of~$Y(\\)$ is equal to
$\dim(\ker\eta^{j+1})-\dim(\ker\eta^j)$ by proposition~\Jordanprop, so from
(\lastlabel) one gets $J(\eta|_{V'})\tr_j\leq\\\tr_j$, and combining
this for all~$j$ yields $J(\eta|_{V'})\subset\\$. Similarly, the kernel
of the projection $\im\eta^j\to(\im\eta^j)/V'=\im(\eta\[V'])^j$ is equal to
$V'\thru\im\eta^j$, so its dimension $\dim(\im\eta^j)-\dim(\im(\eta\[V'])^j)$
decreases weakly as $j$ increases, since $\im\eta^j\supset\im\eta^{j+1}$. Then
using $\\\tr_j=\dim(\im\eta^j)-\dim(\im\eta^{j+1})$, it follows analogously to
the argument above that $\\\tr_j\geq J(\eta\[V'])\tr_j$ for all~$j$,
whence $J(\eta\[V'])\subset\\$.
\QED

The partition $J(\eta|_{V'})$ is called the type of~$V'$, and $J(\eta\[V'])$
is the cotype of~$V'$ (in~$V$). Since the spaces
$\ker(\eta|_{V'})^j=V'\thru\ker\eta^j$ whose dimensions determine the type
of~$V'$ are not directly related to the spaces $V'\thru\im\eta^j$ that were
used to determine its cotype, it is not generally possible (for a fixed value
of~$\\$) to determine the type from the cotype or vice versa. However, there
is an exception when $\\$ is a ``rectangular'' partition, i.e., when all its
non-zero parts have a fixed size~$d$: in that case one has
$\ker\eta^j=\im\eta^{d-j}$ for $0\leq j\leq d$. This leads to the following
fact.

\proclaim Proposition. \rectJordanprop
If $\\=(d,d,\ldots,d)$ with $d$ occurring $m$ times, and
$J(\eta|_{V'})=(\li(\mu0..m-1))$, then
$J(\eta\[V'])=(d-\mu_{m-1},\ldots,d-\mu_0)$.

\proof From the proof of proposition~\Jsubsetprop, one has for $0\leq j<d$:
$$ \displaylines {\quad
 J(\eta|_{V'})\tr_j
 =\dim(V'\thru\ker\eta^{j+1})-\dim(V'\thru\ker\eta^j)
\hfill\cr\hfill{}
 =\dim(V'\thru\im\eta^{d-j-1})-\dim(V'\thru\im\eta^{d-j})
 =\\\tr_{d-j-1}-J(\eta\[V'])\tr_{d-j-1}
 =m-J(\eta\[V'])\tr_{d-j-1},
\quad\cr}
$$
from which the stated relation between $J(\eta|_{V'})$ and $J(\eta\[V'])$
follows.
\QED

We now specialise to the cases of $\eta$-stable lines and hyperplanes. As we
have seen above, $\eta$-stability of a line $l$ means $l\subset\ker\eta$, so
the set of $\eta$-stable lines is identified with the projective space
$\P(\ker\eta)$, which we shall denote by~$\P(V)_\eta$. Also, the cotype
$J(\eta\[l])$ of~$l$ is determined by the (weakly decreasing) sequence of
values $\dim(l\thru\im\eta^j)\in\{0,1\}$, for $j\geq0$. We define for $j\in\N$
subspaces
$$
  W_j(\eta)=\ker\eta\thru\im\eta^j \nn
$$
of~$\ker\eta$, which form a weakly decreasing chain. We also define subsets
$$
  U_j(\eta)=\P(W_j(\eta))\setminus\P(W_{j+1}(\eta)) \nn
$$
of the projective space~$\P(V)_\eta$; the non-empty $U_j(\eta)$ form a finite
partition of that space. We have $\dim W_j(\eta)=\\\tr_j$ by
proposition~\Jordanprop; therefore $U_j(\eta)$ is non-empty if and only if the
following equivalent statements hold: $\\\tr_j>\\\tr_{j+1}$; there is a corner
of~$\\$ in column~$j$; at least one part of~$\\$ equals $j+1$. For the case of
$\eta$-stable hyperplanes we can apply these definitions to $\eta^*$, the nilpotent
endomorphism of the dual vector space~$V^*$ given by
$\eta^*(\phi)\:v\mapsto\phi(\eta(v))$ for $\phi\in V^*$ and $v\in V$. We
therefore define
$$
  W^*_j(\eta)=W_j(\eta^*)
 =\setof\phi\in V^*: \phi(\im\eta)=0 \land \phi(\ker\eta^j)=0\endset. \nn
$$
The set $U_j(\eta^*)$ is contained in the set~$\P(V^*)$ of $1$-dimensional
subspaces of~$V^*$, which is in canonical bijection with the set of
hyperplanes~$H$ of~$V$ by $H\mapsto\setof\phi\in V^*:\phi(H)=0\endset$. We
shall denote this set of hyperplanes by~$\P^*(V)$, and its subset
$\setof H\in\P^*(V):H\supset\im\eta\endset$ of $\eta$-stable hyperplanes
by~$\P^*(V)_\eta$. Then we define $U^*_j(\eta)$ as the subset
of~$\P^*(V)_\eta$ corresponding to~$U_j(\eta^*)$:
$$ 
  U^*_j(\eta)=
  \setof H\in\P^*(V)_\eta:H\supset\ker\eta^j\land H\not\supset\ker\eta^{j+1}
  \endset. \nn
$$
The $U_j(\eta)$ and $U^*_j(\eta)$ respectively partition $\P(V)_\eta$
according to cotype and $\P^*(V)_\eta$ according to type:

\proclaim Proposition. \boxremoveprop
If\/ $l\in U_j(\eta)$, then the Young diagram of $J(\eta\[l])$ is obtained
from that of~$\\$ by removing its corner in column~$j$. If\/ $H\in
U^*_j(\eta)$, then the Young diagram of $J(\eta|_H)$ is obtained from that
of~$\\$ by removing its corner in column~$j$.

\proof
If $l\in U_j(\eta)$, then by reasoning as in the proof of
proposition~\Jsubsetprop\ we find that $\\\tr_j-J(\eta\[V'])\tr_j=1$, and
$\\\tr_c=J(\eta\[V'])\tr_c$ for all $c\neq j$. The argument for
$H\in U^*_j(\eta)$ is entirely analogous.
\QED

For any basis $\{b_0,\ldots,b_k\}$ of~$\ker\eta$ with the property that each
$W_j(\eta)$ is spanned by $\setof b_i:0\leq i<\\\tr_j\endset$, one can find a
decomposition into Jordan blocks $V=B_0\directsum\cdots\directsum B_k$ such
that $\ker(\eta|_{B_i})=\<b_i>$ for all~$i$ (it~suffices to choose
vectors~$v'_i$ with $\eta^j(v'_i)=b_i$, where $\<b_i>\in U_j(\eta)$, and set
$B_i=\<\eta^k(v'_i)\mid0\leq k\leq j>$). For any given $l\in\P(V)_\eta$ the
basis can be chosen such that $l=\<b_i>$ for some~$i$; we shall call a
corresponding decomposition of~$V$ into Jordan blocks adapted to~$l$. We shall
similarly call a decomposition of~$V$ into Jordan blocks adapted
to~$H\in\P^*(V)_\eta$ if $H$ contains all these blocks but~one; for that
block~$B_i$ one has $H\thru B_i=\im(\eta|_{B_i})$. We see
that the centraliser $Z_u$ of~$u$ in~$\GL(V)$ acts transitively on each
set~$U_j(\eta)$ and on each~$U^*_j(\eta)$. The following characterisations of
the index~$j$ such that $l\in U_j(\eta)$ respectively $H\in U^*_j(\eta)$ will
be useful in the sequel.

\proclaim Lemma. \Hlauxlemma
\statitem If $l\in U_j(\eta)$, then $j$ is the minimal value for which
     $(\ker\eta)/l\supset W_j(\eta\[l])$.
\statitem If $H\in U^*_j(\eta)$, then $j$ is the minimal value for which
     $\im\eta\subset\im(\eta|_H)+\ker(\eta|_H)^j$.
\statitem
     $H\in U^*_j(\eta)$ if and only if
     $\im\eta+\ker\eta^j=\im(\eta|_H)+\ker(\eta|_H)^j$.
\statitem If $H\in U^*_j(\eta)$, then $W_c(\eta|_H)=W_c(\eta)$ for
     $c\neq j$, while $W_j(\eta|_H)$ is a hyperplane in~$W_j(\eta)$.

\def\keretalemma{\Hlauxlemma\thinspace\statitemnr1}
\def\imetalemma{\Hlauxlemma\thinspace\statitemnr2}
\def\imkerlemma{\Hlauxlemma\thinspace\statitemnr3}
\def\Wlemma{\Hlauxlemma\thinspace\statitemnr4}

\proof
In each case let the initial condition be satisfied, and let
$V=B_0\directsum\cdots\directsum B_k$ be a decomposition into Jordan blocks
adapted to~$l$ respectively to~$H$. For~\statitemnr1, let $B_i$ be the block
containing~$l$, then $V/l\isom
B_0\directsum\cdots\directsum(B_i/l)\directsum\cdots\directsum B_k$, and the
projection $V\to V/l$ is the identity on all summands except~$B_i$. This
reduces us to the case $V=B_i$; then $(\ker\eta)/l=\{0\}$ and
$J(\eta\[l])=(j)$, whence the statement is obvious. For~\statitemnr2, let
$B_i$ be the block not contained in~$H$, then the intersection of~$\im\eta$
with any other block is contained in~$\im(\eta|_H)$; this again reduces us to
the case $V=B_i$, where the statement follows from $J(\eta|_H)=(j)$.
Part~\statitemnr3 now follows because $\ker\eta^c\not\subset H$ for~$c>j$.
Part~\statitemnr4 also follows by considering the decomposition of~$V$, or by
observing that $W_c(\eta|_H)\subset W_c(\eta)$ and $\dim(W_c(\eta))=\\\tr_c$,
in conjunction with proposition~\boxremoveprop.
\QED

 \par
\newsection Flags.
\labelsec\flagsec

A (complete) \newterm{flag}~$f$ in~$V$ is a saturated chain $0=f_0\ssubset f_1
\ssubset \cdots\ssubset f_n=V$ of subspaces of~$V$. We have $\dim f_i=i$, and
the individual spaces~$f_i$ are called the \newterm{parts} of~$f$. Let $\F$ be
the set of all such flags, called the \newterm{flag manifold} of~$V$. It has
the structure of a projective algebraic variety (see
\ref{Humphreys,~8.1}), and the maps $f\mapsto f_i$ are morphisms onto the
respective Grassmann varieties. Of particular interest are the morphisms giving
the line and hyperplane parts: if $n>0$ we define $\alpha\:\F\to\P(V)$ and
$\omega\:\F\to\P^*(V)$ by $\alpha(f)=f_1$, $\omega(f)=f_{n-1}$. The group
$\GL(V)$ acts on~$\F$, $\P(V)$, and~$\P^*(V)$, and clearly
$\alpha$~and~$\omega$ are $\GL(V)$-equivariant. We say that a flag $f\in\F$ is
$\eta$-stable if all its parts are, and let $\Fu$ denote the subvariety of
$\eta$-stable flags in~$\F$, or equivalently the fixed point set of~$u$ acting
on~$\F$; $\au$~and~$\ou$ will denote the restrictions to~$\Fu$ of
$\alpha$~and~$\omega$, respectively. We have
$\im\au=\P(V)_\eta=\Union_{j\geq0}U_j(\eta)$, and similarly
$\im\ou=\P^*(V)_\eta=\Union_{j\geq0}U^*_j(\eta)$.

For each $\eta$-stable hyperplane~$H\in\P^*(V)_\eta$, the fibre~$\ou\inv(H)$
of~$\ou$ is isomorphic to the variety~$\FuH$ of $\eta|_H$-stable flags in~$H$.
Indeed the \iso. is given by $f\mapsto f^-$, where
$f^-=(f_0\ssubset\cdots\ssubset f_{n-1}=H)$ is the flag obtained from~$f$ by
omitting the last part $f_n=V$. Similarly, for each $\eta$-stable
line~$l\in\P(V)_\eta$, the fibre~$\au\inv(l)$ is isomorphic to the variety
$\Ful$ of $\eta\[l]$-stable flags in~$V/l$; here the \iso. will be written as
$f\mapsto f"$, where $f"=(f_1/l\ssubset\cdots\ssubset f_n/l)$ is the flag
obtained from~$f$ by reducing modulo~$f_1=l$ all its parts except~$f_0$. There
will be no confusion if the same notations $f^-$ and~$f"$ are used when $f$ is
a flag in another vector space than~$V$ (provided its dimension is non-zero);
this allows us in particular to write for~$f\in\Fu$ expressions such as
$f^{--}$, $f"""$, and $f^{-\downarrow}$, as long as the total number of
operations applied does not exceed~$n$. Moreover, the operations commute:
$f^{-\downarrow}$ and $f"^-$ denote the same flag in~$f_{n-1}/f_1$. In
general, one obtains in this manner from~$f\in\Fu$ an $\eta_{f_i/f_j}$-stable
flag in $f_i/f_j$, for some $i\geq j$, where $\eta_{f_i/f_j}$ is the nilpotent
endomorphism of~$f_i/f_j$ induced by~$\eta$.

The sequences of types and of cotypes of the parts of~$f\in\Fu$ define two
saturated decreasing chains in~$\Part$ from~$\\$ to~$(0)$, that can be used to
define Young tableaux $r_\eta(f),q_\eta(f)\in\Tab_\\$:
$$ \eqalignno{
\ch r_\eta(f)&=\( J(\eta),J(\eta|_{f_{n-1}}),J(\eta|_{f_{n-2}}),\ldots,(0)\)
   & \nn \cr
\ch q_\eta(f)&=\(J(\eta),J(\eta\[f_1]),J(\eta\[f_2]),\ldots,(0)\).
   & \nn \cr}
$$
In other words, the subtableau of $r_\eta(f)$ containing entries~$<i$ has
shape~$J(u|_{f_i})$, while the subtableau of $q_\eta(f)$ containing
entries~$<n-i$ has shape~$J(u\[f_i])$. If we define for each flag~$f\in\F$ a
dual flag~$f^*$ in~$V^*$ by ${f^*}_i=\setof\phi\in V^*:
\phi(f_{n-i})=0\endset$, then one readily verifies that
$q_\eta(f)=r_{\eta^*}(f^*)$, and $r_\eta(f)=q_{\eta^*}(f^*)$.

As was mentioned earlier, there is in general no direct relationship between
the type and cotype of the parts of~$f$, and so there is no one-to-one
correspondence between $r_\eta(f)$ and~$q_\eta(f)$ either. However, we have
again an exception if $\\$ is a rectangular partition. To describe the
relationship in this case, we introduce the involutive operation $T\mapsto
T^\diamond$ on~$\Tab_\\$ for rectangular~$\\$: let the square~$t$ be the
unique corner of~$\\$, then whenever some square $s\leq t$ has entry~$i$
in~$t$, then the diametrically opposite square~$t-s$ has entry $\w_i=n-1-i$
in~$T^\diamond$ (this is essentially the same as the operation
$P\mapsto\overline P$ of \ref{festschrift, proposition~5.7}).

\proclaim Proposition. \qrThm
If $\\$ is a rectangular partition, then $q_\eta(f)=r_\eta(f)^\diamond$ for
all $f\in\Fu$.

\proof
Put $\ch r_\eta(f)= (\\^n,\\^{n-1},\ldots,\\^0)$, and $\ch q_\eta(f)=
(\mu^n,\mu^{n-1},\ldots,\mu^0)$. Then each $\\^i$ determines $\mu^{n-i}$, as
described in proposition~\rectJordanprop.  For $0\leq i<n$, the square with
entry~$i$ in~$r_\eta(f)$ is $\\^{i+1}-\\^i$, which determines the square
$\mu^{n-i}-\mu^{n-i-1}$ containing $\w_i$ in~$q_\eta(f)$; therefore
$q_\eta(f)=r_\eta(f)^\diamond$.
\QED

We define for any Young tableau~$T$ of shape~$\\$:
$$ \eqalignno
{ \F_{\eta,T}  &=\setof f\in\Fu: r_\eta(f)\sim T \endset & \nn \cr
  \F^*_{\eta,T}&=\setof f\in\Fu: q_\eta(f)\sim T \endset. & \nn \cr
}
$$
If the square~$\hi T$ lies in column~$j$ (which implies that $U_j(\eta)$ and
$U^*_j(\eta)$ are non-empty), then one has $\omega(\F_{\eta,T})\subset
U^*_j(\eta)$ and $\alpha(\F^*_{\eta,T})\subset U_j(\eta)$. Moreover, for any
$H\in U^*_j(\eta)$ the fibre $\F_{\eta,T}\thru\ou\inv(H)$ is isomorphic to
$\F_{\eta|_H,T^-}$ by $f\mapsto f^-$, and for any $l\in U_j(\eta)$ the fibre
$\F^*_{\eta,T}\thru\au\inv(l)$ is isomorphic to $\F^*_{\eta\[l],T^-}$ by
$f\mapsto f"$. It follows by induction that each of the sets $\F_{\eta,T}$ and
$\F^*_{\eta,T}$ is non-empty, and open in its closure. As $T$ ranges
over~$\Tab_\\$, the sets $\F_{\eta,T}$ partition~$\Fu$ into finitely many
subsets, as do the sets~$\F^*_{\eta,T}$.

\proclaim Proposition. \FuTirrprop
\statitem
  For each $T\in\Tab_\\$ the sets $\F_{\eta,T}$ and $\F^*_{\eta,T}$ are
  irreducible.
\statitem
  $\dim(\F_{\eta,T})=\dim(\F^*_{\eta,T})=n(\\)\defeq\sum_{i\geq0}i\\_i$,
  independently of\/~$T\in\Tab_\\$.
\statitem
  The set of irreducible components of~$\Fu$ is equal to
  $\setof\overline{\F_{\eta,T}}:T\in\Tab_\\\endset$, and also to
  $\setof\overline{\F^*_{\eta,T}}:T\in\Tab_\\\endset$.

\proof
We proceed by induction on $n=|\\|$, and only prove the statements
for~$\F_{\eta,T}$, as those for $\F^*_{\eta,T}$ are entirely similar (and also
follow by transition to the dual vector space). The case $n=0$ is trivial, so
assume $n>0$; let $T\in\Tab_\\$ and let $\hi T$ be the square~$(i,j)$. Since
$U^*_j(\eta)$ is irreducible and an orbit for~$Z_u$, it is already an orbit
for the identity component $Z_u\o$ of~$Z_u$ (in fact $Z_u$ is always
connected, but we do not wish to invoke that fact here). Using the \iso.
$\F_{\eta,T}\thru\ou\inv(H)\isomto\F_{\eta|_H,T^-}$ for some $H\in
U^*_j(\eta)$, we may define a surjective morphism
$Z_u\o\times\F_{\eta|_H,T^-}\to\F_{\eta,T}$ by $(z,f^-)\mapsto z\cdot f$;
since the domain of this morphism is irreducible by the induction hypothesis,
so is its image, which establishes~\statitemnr1. We have
$\dim(U^*_j(\eta))=\dim(W^*_j(\eta))-1=i$, and so for any $H\in U^*_j(\eta)$
we have $\dim\F_{\eta,T}=\dim U^*_j(\eta)+\dim\F_{\eta|_H,T^-}=i+n(\sh
T^-)=n(\\)$, proving~\statitemnr2. Part~\statitemnr3 follows from
\statitemnr1~and~\statitemnr2.
\QED

\heading "Remark"
Part~\statitemnr3 gives two different natural parametrisations of the
irreducible components of~$\Fu$ by Young tableaux. The first one, based on the
types of the parts of flags (as it uses~$r_\eta$) corresponds to the
parametrisation used in~\ref{SteinbergORS}, but in~\ref{Spaltenstein,~II.5.3}
the other parametrisation, based on cotypes ($q_\eta$), is effectively used.
We choose to work primarily with the former parametrisation, partly because it
is somewhat simpler to use restrictions than quotients, but mainly because
this choice leads to a more direct interpretation of the Robinson-Schensted
algorithm. Note that our choice does lead to a slightly illogical use of
asterisks: $\F_{\eta,T}$ has a fibration over $U^*_j(\eta)$, while
$\F^*_{\eta,T}$ has one over~$U_j(\eta)$.

We close this section with an example,
illustrating these parametrisations of the irreducible components of~$\Fu$
in the simplest non-trivial case, namely for the Jordan type~$\\=(2,1)$.
Then $\Tab_\\$ has $2$~elements, namely
$$
   T=\Young(0,1|2) \ttext{and} T'=\Young(0,2|1),
$$
and hence $\Fu$ has $2$~irreducible components, of dimension $n(\\)=1$. To be
specific, let us take
$$
  \eta=\pmatrix{ 0 & 1 & 0 \cr 0 & 0 & 0 \cr 0 & 0 & 0 \cr}.
$$
Calling the standard basis vectors $e_0$,~$e_1$,~$e_2$ we have
$\ker\eta=\<e_0,e_2>$ and $\im\eta=\<e_0>$, while $\eta^j=0$ for $j>1$. There
are two $Z_u$-orbits of $\eta$-stable hyperplanes, namely the set $U^*_0(\eta)$
of all hyperplanes containing $\im\eta=\<e_0>$ but not $\ker\eta=\<e_0,e_2>$,
and the singleton $U^*_1(\eta)=\{\<e_0,e_2>\}$. For any $H\in U^*_0(\eta)$ the
fibre~$\ou\inv(H)$ consists of a single flag~$f$, with
$f_1=H\thru\ker\eta=\<e_0>$ and of course $f_2=H$, so $\F_{\eta,T}$ is
isomorphic to the affine line~$U^*_0(\eta)$. On the other hand for
$H=\<e_0,e_2>\in U^*_1(\eta)$ we have $\eta|_H=0$, whence an $\eta$-stable
flag~$f$ with $f_2=H$ can have an arbitrary element of~$\P(H)$ as $f_1$, so in
this case the fibre~$\ou\inv(H)$, which coincides with $\F_{\eta,T'}$, is a
projective line. There is one element of~$\F_{\eta,T'}$ that lies in the
closure of $\F_{\eta,T}$, namely the flag~$f$ with $f_1=\<e_0>$ and
$f_2=\<e_0,e_2>$. Therefore, the whole variety~$\Fu$ can be depicted as
follows
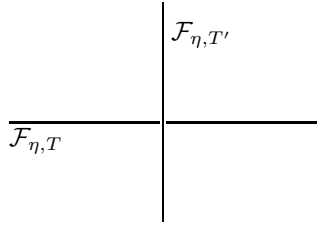
\bigdisplay
\vcenter{\offinterlineskip\halign{& # \hfil&\kern1pt\vrule#\kern1pt\cr
  & height 5mm
  & \kern2pt$\F_{\eta,T'}$ \cr
  \vrule height.4pt depth0pt width 2cm & height 1cm depth 2pt &
  \vrule height.4pt depth0pt width 2cm \cr
  $\F_{\eta,T}$ & depth 1cm \cr}}
$$
This illustration should be considered to lie in $\P_1\times\P_1$, with the
horizontal coordinate representing the choice of the hyperplane $H=f_2$, and
the vertical coordinate representing the choice of the line~$f_1$. By similar
reasoning to that above, it can be seen that the fibre $\F^*_{\eta,T'}$ of the
morphism $\au$ is the horizontal line \emph{including} the point of crossing,
and $\F^*_{\eta,T}$ is the remainder of the vertical line. So $r_\eta(f)\neq
q_\eta(f)$ for all flags $f\in\Fu$, except the flag represented by the
crossing point of the lines, for which
$r_\eta(f)=q_\eta(f)=T'={\smallsquares\smallsquares\Young(0,2|1)}$.
 \par
\newsection {Interpretation of the \Sch. algorithm}.

We shall now proceed to show that the two given parametrisations of the
irreducible components of~$\Fu$ are related by the involution~$S$ of~$\Tab_\\$
defined by the \Sch. (evacuation) algorithm. In fact we shall give a more
detailed description of the situation than this. We first recall the following
two kinds of configurations of $4$~partitions that can occur within the doubly
indexed family of partitions that was used in \ref{festschrift,~2.2} to give a
description of the \Sch. algorithm.

\proclaim\De. \Sconfigdef
An arrangement of 4~partitions $\ar(\\,\mu,\mu',\nu)$ with
$\nu\prec\mu,\mu'\prec\\$ is called
\item{(i)}  a configuration of type~$S1$ if $\mu=\mu'$ and
            $(\mu-\nu)\adj(\\-\mu)$,
\item{(ii)} a configuration of type~$S2$ if $\mu\neq\mu'$ (so that one has
            $\\=\mu\union\mu'$ and $\nu=\mu\thru\mu'$).

Note that the condition $(\mu-\nu)\adj(\\-\mu)$ implies $\mu=\mu'$ (there are
no other partitions between $\nu$~and~$\\$), and that in the other case the
two partitions $\mu,\mu'$ are the only ones between $\nu$~and~$\\$. It follows
that if $\ar(\\,\mu,\mu',\nu)$ is known to be of type $S1$~or~$S2$, then
either one of $\mu$~or~$\mu'$ is uniquely determined by the other three
partitions. The relevance of these configurations becomes clear when one
considers for $\nu\prec\mu\prec\\$ the following variety $\E_{\eta,\mu,\nu}$
of partial flags, consisting of only a line and hyperplane part:
$$
  \E_{\eta,\mu,\nu}=\setof (l,H)\in \P(V)_\eta\times\P^*(V)_\eta:
   J(\eta|_H)=\mu \land l\subset H \land J(\eta_{H/l})=\nu \endset.
$$

\proclaim Lemma. \Elemma
Let $\nu\prec\mu\prec\\$, and let $\mu'$ be determined by the
condition that $\ar(\\,\mu,\mu',\nu)$ is of type $S1$~or~$S2$. Then
$\E_{\eta,\mu,\nu}$ is irreducible, and $J(\eta\[l])=\mu'$ for $(l,H)$ in a
$Z_u$-stable dense open subset of $\E_{\eta,\mu,\nu}$.

\proof
Put $\\-\mu=(i,j)$, $\mu-\nu=(r,c)$; then $\E_{\eta,\mu,\nu}$ consists of
those pairs $(l,H)$ with $H\in U^*_j(\eta)$ and $l\in U_c(\eta|_H)$, from
which its irreducibility follows. If $j=c$ or $i=r$, then
$\ar(\\,\mu,\mu',\nu)$ is of type~$S1$, and $J(\eta\[l])=\mu'$ is forced by
$\nu\prec J(\eta\[l])\prec\\$. If $j\neq c$ and $i\neq r$, then
$\ar(\\,\mu,\mu',\nu)$ is of type~$S2$; if moreover $c\neq j-1$, then
$U_c(\eta|_H)=U_c(\eta)$ follows from lemma~\Wlemma, so
$J(\eta\[l])=\\-(r,c)=\mu'$. In the final case that $c=j-1$ and $r>i$, we find
for all $H\in U^*_j(\eta)$, again using lemma~\Wlemma, that
$U_c(\eta)=\P(W_c(\eta))\setminus\P(W_j(\eta))\isom\P_r\setminus\P_i$ is dense
and open in $U_c(\eta|_H)=\P(W_c(\eta))\setminus\P(W_j(\eta|_H))\isom
\P_r\setminus\P_{i-1}$ (with $\P_{-1}=\emptyset)$), so the subset
of~$\E_{\eta,\mu,\nu}$ of pairs $(l,H)$ for which $l\in U_c(\eta)$ is dense
and open, and on this set $J(\eta\[l])=\\-(r,c)=\mu'$ holds.
\QED

\proclaim Theorem. \Schinterp
There is a dense open subset $\Fu'$ of~$\Fu$, such that for all $f\in\Fu'$ the
following holds: if for $i+j\leq n$ one puts
$\\^{[i,j]}=J(\eta_{f_{n-j}/f_i})$, then any configuration
$$\Ar(\\^{[i,j]},\\^{[i,j+1]},\\^{[i+1,j]},\\^{[i+1,j+1]})$$
is either of type~$S1$ or of type~$S2$.

\proof
It will suffice to find for each $T\in\Tab_\\$ a dense open subset
of~$\F_{\eta,T}$ on which the stated condition holds. We use induction on~$n$,
with $n=1$ as trivial starting case. Fix $T\in\Tab_\\$, and let the family
$(\\^{[i,j]})_{i+j\leq n}$ of partitions be defined by $\ch T=
(\\^{[0,0]},\\^{[0,1]},\ldots,\\^{[0,n]})$, and by the condition on
configurations in the statement of the theorem; we shall show that for
appropriately chosen $f\in\F_{\eta,T}$ one has
$J(\eta_{f_{n-j}/f_i})=\\^{[i,j]}$ for $i+j\leq n$. Choose some
$H\in\omega(\F_{\eta,T})$, then by the induction hypothesis applied to
$\F_{\eta|_H,T^-}$, we have for $f$ in a dense open subset of the fibre
$\F_{\eta,T}\thru\ou\inv(H)$ that all instances of
$J(\eta_{f_{n-j}/f_i})=\\^{[i,j]}$ with~$j>0$ hold. Since
$\omega(\F_{\eta,T})$ is equal to some $Z_u$-orbit $U^*_c(\eta)$, this
property extends by the action of~$Z_u$ to a $Z_u$-stable dense open
subset~$\U$ of~$\F_{\eta,T}$. By construction $\U$ is contained in $\setof
f\in\F_{\eta,T}:(f_1,f_{n-1})\in\E_{\eta,\mu,\nu}\endset$, which projects onto
$\E_{\eta,\mu,\nu}$ with fibres that are all of dimension~$n(\\^{[1,1]})$. It
follows that $\U$ meets the inverse image under this projection of the dense
open subset of~$\E_{\eta,\mu,\nu}$ described in lemma~\Elemma. The
intersection~$\U'$ of~$\U$ with this inverse image is a $Z_u$-stable dense
open subset of~$\F_{\eta,T}$ on which the additional property
$J(\eta\[f_1])=\\^{[1,0]}$ holds. Then $\alpha(\U')$ is a single $Z_u$-orbit,
and $\U'$ intersects each fibre $\au\inv(l)$ for $l\in\alpha(\U')$ densely;
the image of this intersection under $f\mapsto f"$ is a dense open subset
of~$\F_{\eta\[l],T"}$. The proof can now be completed by applying the
induction hypothesis to~$\F_{\eta\[l],T"}$.
\QED

\proclaim Corollary. \paramchange
For all $T\in\Tab_\\$ one has
$\overline{\F_{\eta,T}}=\overline{\F^*_{\eta,S(T)}}$.

\proof
This follows immediately from the theorem, and the fact the the family
$\\^{[i,j]}$ with the mentioned properties can be used to compute $S(T)$
(cf.~the proof of \ref{festschrift, theorem~2.2.1}).
\QED

The proof given for this geometric interpretation uses a description of~$S$
from which it obvious that $S$ is an involution (in fact one that was used to
prove this fact combinatorially). However, even if we forget the proof, the
interpretation clearly implies that $S$ is an involution, since
$f\in\F^*_{\eta,T}$ is equivalent to $f^*\in\F_{\eta^*,T}$, and $f^{**}=f$. In
combination with proposition~\qrThm, the corollary also implies that $S(T)\sim
T^\diamond$ for tableaux~$T$ of rectangular shape, as stated
in~\ref{festschrift, Corollary~5.7}; in this case the geometric proof has no
purely combinatorial counterpart (the proof given in~\ref{festschrift} is
based on the relationship of the \Sch. correspondence with the
Robinson-Schensted correspondence). There is on the other hand one
combinatorially obvious symmetry of~$S$ without any clear geometric meaning,
namely the fact that it commutes with transposition (a geometric
interpretation of this fact would require some operation that gives rise to
transposition of Jordan types). Nevertheless, this symmetry of~$S$ will be
important below, when we give geometric interpretations of the
Robinson-Schensted correspondence and its transpose.

 \par
\newsection
Interpretation of jeu de taquin, and Littlewood-Richardson tableaux.

The deflation procedure used in the \Sch. algorithm is related to the
operation of jeu de taquin (also known as \foreign{glissement}), which is
performed on skew tableaux. Jeu de taquin can be described completely in terms
of the deflation procedure \ref{festschrift, \Sec5}, which allows us to deduce
from theorem~\Schinterp\ an interpretation of jeu de taquin. That theorem does
not mention tableaux or their entries directly however, but is stated in terms
of families of partitions, and for our interpretation, all that matters about
a skew tableau is the chain in~$\Part$ associated to it. For convenience we
shall work directly with such chains, rather than with skew tableaux: we
define a skew chain of shape $\\/\mu$ to be a saturated decreasing chain
in~$\Part$ from~$\\$ to~$\mu$. We denote the set of all skew chains of shape
$\\/\mu$ by $\Tab_{\\/\mu}$; the set $\Tab_\\$ is in bijection with
$\Tab_{\\/\emptyset}$ by $P\mapsto\ch P$.

We proceed to give a geometric interpretation to skew chains, in analogy to
the definition of~$\F_{\eta,T}$ for Young tableaux~$T$. Let $\Gr{m}$ be the
variety of $m$-dimensional $\eta$-stable subspaces of~$V$, and let
$\F^{(m)}_\eta$ be the set of $\eta$-stable partial flags in~$V$ with parts in
dimensions $m$ and higher, i.e., of chains $f={(f_m\ssubset
f_{m+1}\ssubset\cdots\ssubset f_n=V)}$ with $f_i\in\Gr{i}$. For
$f\in\F^{(m)}_\eta$ the part~$f_m$ of minimal dimension will be denoted
by~$\lo f$, and we define the complete flag $\overline f=(f_m/\lo
f\ssubset\cdots\ssubset f_n/\lo f)\in\F_{\eta\[\lo f]}$ by reducing all parts
of~$f$ modulo~$\lo f$. We also define $r_\eta(f)=
\(J(\eta),J(\eta|_{f_{n-1}}),\ldots,J(\eta|_{\lo f})\)\in\Tab_{\\/\mu}$ where
$\mu=J(\eta|_{\lo f})$, and put
$$
  \F_{\eta,K}= \setof f\in\F^{(m)}_\eta: r_\eta(f)=K\endset.
$$

\proclaim Proposition.
$\F_{\eta,K}$ is an irreducible variety of dimension $n(\\)-n(\mu)$, for every
$K\in\Tab_{\\/\mu}$.

\proof
The proof is entirely analogous to that of parts (a)~and~(b) of
proposition~\FuTirrprop, the only difference being that the induction starts
at $|\\|=|\mu|$ rather than at $|\\|=0$.
\QED

For a Young tableau $T\in\Tab_\\$, let $T_{<m}$ denote its subtableau of
entries less than~$m$; putting $\mu=\sh T_{<m}$, let
$T_{\geq m}\in\Tab_{\\/\mu}$ denote the subchain $(\sh T, \sh T^-,\ldots,\mu)$
of $\ch T$ corresponding to the remaining entries. If we denote
by~$T^{\downarrow*i}$ the result of applying the deflation procedure $i$~times
to~$T$, then the relation $K\slid K'$ of glissement ($K'$~is obtainable
from~$K$ by inward jeu de taquin slides) is defined in \ref{festschrift,
\Sec5} by restricting $T$ and $T^{\downarrow*i}$ to skew subtableaux with the
same set of entries. The corresponding definition for skew chains is that
$T_{\geq{m}}\slid{T^{\downarrow*i}}_{\geq m}$ for $i\leq m$. We obtain from
theorem~\Schinterp:

\proclaim Corollary. \jdtinterp
Let $K\in\Tab_{\\/\mu}$ be a skew chain, and let $P\in\Tab_\nu$ be a Young tableau
such that $K\slid\ch{P}$. Then there is a dense open subset $\F'_{\eta,K}$ of
$\F_{\eta,K}$ such that $r_{\eta\[\lo f]}(\overline f)=P$ for
all~$f\in\F'_{\eta,K}$.

\proof
Since $K\slid\ch P$, there exists a $T\in\Tab_\\$ with $K=T_{\geq m}$ and
$P\sim T^{\downarrow*m}$. Let $p$ denote the natural projection
$\F_{\eta,T}\to\F_{\eta,K}$ (which is clearly surjective), and take
$\F'_{\eta,K}=p(\F_{\eta,T}\thru\Fu')$, where $\Fu'$ is as in
theorem~\Schinterp; the conclusion of that theorem immediately implies
$r_{\eta\[\lo f]}(\overline f)=P$ for~$f\in\F'_{\eta,K}$.
\QED

In this geometric interpretation, confluence of jeu de taquin is obvious: $P$
is uniquely determined by~$K$. We now partition $\Gr m$ according to the type
and cotype of its elements: for $\mu\in\Part_m$, $\nu\in\Part_{n-m}$ put
$$
  \Gx=\setof X\in\Gr{m}: J(\eta|_X)=\mu \land J(\eta\[X])=\nu \endset;
$$
furthermore, denote by~$\Tx$ the set of $K\in\Tab_{\\/\mu}$ such that
$K\slid\ch P$ for some~$P\in\Tab_\nu$.

\proclaim Theorem. \LRinterp
Denote by $\pi\:\F^{(m)}_\eta\to\Gr{m}$ the morphism given by $f\mapsto\lo f$.
Let $\mu\in\Part_m$ and $\nu\in\Part_{n-m}$.
\statitem
  If $K\in\Tx$, then the image $\pi(\F'_{\eta,K})$ of the set $\F'_{\eta,K}$
  of corollary~\jdtinterp\ is dense in an irreducible component of $\Gx$ of
  dimension~$n(\\)-n(\mu)-n(\nu)$. Moreover, any such component is so obtained
  as $\pi(\overline{\F_{\eta,K}})\thru\Gx$ for at least one~$K\in\Tx$, and
  there are no components of higher dimension.
\statitem
  The number of irreducible components of $\Gx$ of
  dimension~$n(\\)-n(\mu)-n(\nu)$ is equal to the Littlewood-Richardson
  coefficient~$c^\\_{\mu,\nu}$, and they can be explicitly parametrised by the
  Littlewood-Richardson tableaux of shape $\\/\mu$ and weight~$\nu$.
\statitem
  For $K,L\in\Tx$, the irreducible components
  $\pi(\overline{\F_{\eta,K}})\thru\Gx$ and
  $\pi(\overline{\F_{\eta,L}})\thru\Gx$ of $\Gx$ are equal if and only if
  $K$~and~$L$ are dual equivalent in the sense of~\ref{Haiman dual
  equivalence}.

Limiting itself to what can be deduced from corollary~\jdtinterp, the theorem
avoids any statement about possible irreducible components of~$\Gx$ of dimension
less than $n(\\)-n(\mu)-n(\nu)$. However, we shall show below that such
components do not exist, and so an accordingly simplified and strengthened
form of the theorem does in fact hold.

\proof
\statitemnr1
It is clear that $\pi\inv(\Gx)\subset\Union_{K\in\Tab_{\\/\mu}}\F_{\eta,K}$,
whose components have dimension $n(\\)-n(\mu)$, and the fibre $\pi\inv(X)$ at
any $X\in\Gx$ is isomorphic to $\F_{\eta\[X]}$, whose components have
dimension~$n(\nu)$. Therefore any irreducible component~$C$ of~$\Gx$ can have
dimension at most $n(\\)-n(\mu)-n(\nu)$, and when it has this dimension,
$\pi\inv(C)$ is dense in some union of sets $\F_{\eta,K}$ with
$K\in\Tab_{\\/\mu}$. By corollary~\jdtinterp\ one has
$\pi(\F'_{\eta,K})\subset\Gx$ whenever $K\in\Tx$, from which the claims
follow.

\statitemnr2
Corollary~\jdtinterp\ also implies that $\F_{\eta,K}$ meets any fibre
$\pi\inv(X)$ with $X\in\pi(\F'_{\eta,K})$ in the irreducible component of that
fibre that corresponds to $\F_{\eta\[X],P}$; hence, fixing an arbitrary
$P\in\Tab_\nu$, the irreducible components of $\Gx$ of dimension
$n(\\)-n(\mu)-n(\nu)$ correspond bijectively to the skew chains $K\in\Tx$ with
$K\slid\ch P$. The number of such~$K$ is known to be independent of the choice
of~$P$ (\ref{Schutzenberger CRGS,~(3.7)}) and equal to $c^\\_{\mu,\nu}$
(\ref{Schutzenberger CRGS,~(4.7)}, see also \ref{pictures, Theorem~5.2.5}). In
fact there is a specific $P\in\Tab_\nu$ for which the Littlewood-Richardson
tableaux~$T$ of shape~$\\/\mu$ and weight~$\nu$ correspond directly to the
skew chains $K\in\Tx$ with $K\slid\ch P$. To associate to a semistandard skew
tableau~$T$ a skew chain~$K$, one uses the well known process of
standardisation: to form the chain~$K$ of partitions starting from~$\\$, the
squares of~$T$ are removed by decreasing entries, and among squares with equal
entries from right to left. Jeu de taquin is defined for semistandard tableaux
in such a way that it commutes with this standardisation, and it preserves the
property of being a Littlewood-Richardson tableau. The indicated special
tableau~$P\in\Tab_\nu$ is such that $\ch P$ is the standardisation of the
tableau~$\CT\nu$ of shape~$\nu$ in which each row~$i$ filled with entries~$i$
($\CT\nu$ is the unique Littlewood-Richardson tableau of shape~$\nu$, and it
has weight~$\nu$). Then $K\in\Tab_{\\/\mu}$ is the standardisation of a
Littlewood-Richardson tableau~$T$ of weight~$\nu$ if and only if
$K\slid\ch P$.

\statitemnr3
Define for $K,L\in\Tx$ the equivalence relation $K\cong L$ to mean
$\pi(\overline{\F_{\eta,K}})\thru\Gx=\pi(\overline{\F_{\eta,L}})\thru\Gx$. We
have established above that for any $P\in\Tab_\nu$ the jeu de taquin
equivalence class $\setof K\in\Tx:K\slid\ch P\endset$ is a set of
representatives for the classes for~`$\cong$'. As the members of such a jeu de
taquin equivalence class are mutually dual inequivalent (this is the easy part
of \ref{Haiman dual equivalence, theorem~2.13}), it will suffice to prove that
$K\cong L$ implies the dual equivalence of $K$~and~$L$. We shall establish
this by finding a sequence of jeu de taquin slides that transforms $K$~and~$L$
respectively into $K',L'\in\Tab_{\nu/\emptyset}$, preserving equality of
shapes at each step; being (chains of) Young tableaux of the same shape,
$K'$~and~$L'$ are dual equivalent \ref{Haiman dual equivalence,
corollary~2.5}, which implies dual equivalence of $K$~and~$L$. Since $K\cong
L$, it is possible to choose partial flags $f\in\F'_{\eta,K}$ and
$f'\in\F'_{\eta,L}$ with $\lo{f}=\lo{f'}$; by extending $f$~and~$f'$
identically by suitably chosen parts in dimensions less than~$m$, one obtains
flags $\hat f,\hat f'\in\F'_\eta$ such that $T=r_\eta(\hat f)$ and
$U=r_\eta(\hat f')$ satisfy $T_{\geq m}=K$, \ $U_{\geq m}=L$, and
$T_{<m}=U_{<m}$. Then for $i\leq m$ the shapes of the skew chains
${T^{\downarrow*i}}_{\geq m}$ and~${U^{\downarrow*i}}_{\geq m}$ are both equal
to $J(\eta\[\hat f_i])/J(\eta_{\lo f/\hat f_i})$ (since $\hat f_i=\hat f'_i$),
which gives the required sequence of slides transforming $K$~and~$L$ into
$K'=\ch T^{\downarrow*m}$ respectively into $L'=\ch U^{\downarrow*m}$.
\QED

The detailed statement of the theorem appears to be new. However, the relation
between~$\Gx$ and Littlewood-Richardson coefficients was already indicated in
\ref{Springer Weyl representations, Theorem~4.4}. The setting there is in fact
more general, with a semi-simple linear algebraic group~$G$ replacing~$\GL_n$;
correspondingly, $\Gx$ is replaced by a variety ${\cal X}'_{A,B}(P)$ of
parabolic subgroups, and Littlewood-Richardson coefficients by decomposition
multiplicities for representations induced from the Weyl group~$W'$ of a Levi
factor of~$P$ to the Weyl group~$W$ of~$G$. (Our theorem corresponds only to
maximal parabolic~$P$, but can be extended easily so as to correspond to
arbitrary parabolic subgroups.) The number of irreducible components of~$\Gx$
of dimension $n(\\)-n(\mu)-n(\nu)$ appears in a somewhat disguised form, as a
decomposition multiplicity $n_{A,B,\phi,\psi}$ for a permutation action of a
group $C_G(A)\times C_G(B)$ on the set of irreducible components of~${\cal
X}'_{A,B}(P)$ of that dimension; for $\GL_n$, the group $C_G(A)\times C_G(B)$
is always trivial, and $n_{A,B,\phi,\psi}$ reduces to the number of components
acted (trivially) upon.

Another, rather weaker, connection between the Littlewood-Richardson rule and
the geometry of~$\Fu$ was indicated in~\ref{Srinivasan}. Its main
theorem~(4.2) corresponds to our corollary~\jdtinterp, but it is stated (and
proved) in a somewhat roundabout fashion in terms of permutations, whose link
to geometry is formed by Steinberg's interpretation of the Robinson-Schensted
correspondence (which will be discussed below). In fact, that theorem itself
involves no geometry at all, and it can be proved in a purely combinatorial
manner. Since only a fixed maximal parabolic subgroup~$P$ is considered, no
connection with the geometry of~$\Gx$ is indicated; the Littlewood-Richardson
coefficients, which arise in relation to jeu de taquin in the same way as
above, are only given their traditional representation theoretic
interpretation.

As we have seen, dual equivalence classes in $\Tx$ are in bijection with
Littlewood-Richardson tableaux of shape~$\\/\mu$ and weight~$\nu$. Once the
latter have all been determined, one can construct the set $\setof
K\in\Tx:K\slid\ch P\endset$ effectively, not just for the special tableau~$P$
indicated in the proof above, but for any given $P\in\Tab_\nu$. Such a
construction is given in the proof of \ref{pictures, Theorem~5.2.5} in terms
of the Robinson-Schensted correspondence for ``pictures''; we shall formulate
it here without using pictures. One associates to any skew chain
$K\in\Tab_{\\/\mu}$ a permutation $w(K)$ by concatenating the rows of the skew
tableau corresponding to~$K$, taking them in order from bottom to top. Call
the two Young tableaux $(P,Q)=RS(w(K))$ the $P$-symbol and $Q$-symbol of~$K$;
then the $P$-symbol of~$K$ characterises its jeu de taquin equivalence class
(indeed $K\slid\ch P$), and the $Q$-symbol its dual equivalence class.
Whenever $c^\\_{\mu,\nu}>0$, all tableaux in $\Tab_\nu$ occur as $P$-symbol of
some $K\in\Tx$, but not necessarily as $Q$-symbol. The set of tableaux that do
so occur, is precisely the set $Q(\\/\mu,\nu)$ of $Q$-symbols of
Littlewood-Richardson tableaux of shape~$\\/\mu$ and weight~$\nu$ (where the
$Q$-symbol of a semistandard skew tableau is defined either as the $Q$-symbol
of its standardisation, or directly by concatenating its rows and applying the
version of the Schensted algorithm that allows repeated entries; either way
the $Q$-symbol is a standard tableau). One then has
$$
  \setof w(K):K\in\Tx\land K\slid\ch P\endset=
  \setof RS\inv(P,Q):Q\in Q(\\/\mu,\nu)\endset,
$$
from which the desired set of skew chains~$K$ is readily reconstructed.

Now as promised we shall rule out the possibility that $\Gx$ could have
irreducible components of dimension less than $n(\\)-n(\mu)-n(\nu)$, which
implies in particular that $\Gx=\emptyset$ whenever $c^\\_{\mu,\nu}=0$. This
requires an algebraic construction that associates a Littlewood-Richardson
tableau to any individual element $X\in\Gx$. We essentially use the
construction described in \ref{Macdonald, II~3}, but since our context is dual
to the one considered there, we shall present an adapted version of the
construction and proof.

\proclaim Proposition.
For any $\\,\mu,\nu\in\Part$, the irreducible components of~$\Gx$ are
precisely those described in theorem~\LRinterp, i.e., $\Gx$ has no irreducible
components of dimension less than $n(\\)-n(\mu)-n(\nu)$.

\proof
We shall construct for any $X\in\Gx$ a tableau $K\in\Tx$ with
$X\in\pi(\F_{\eta,K})$; then $X$ lies in the component
$\pi(\overline{\F_{\eta,K}})\thru\Gx$, and the proposition follows. Fix
$X\in\Gx$, and for $i\in\N$ put $X_i=\eta^{-i}(X)$ and $\mu^i=J(\eta|_{X_i})$
(in particular $\mu^0=\mu$ and $\mu^{\nu_0}=\\$). By filling each skew diagram
$Y(\mu^{i+1})\setminus Y(\mu^i)$ with entries~$i$ we obtain the transpose of a
tableau~$T$ of shape~$\\\tr/\mu\tr$ and weight~$\nu\tr$ (by
proposition~\Jsubsetprop, since~$J(\eta\[X])=\nu$); we claim that $T$ is a
Littlewood-Richardson tableau (i.e., $T\slid\CT{\nu\tr}$). Assuming this for
the moment, the transposes $K\in\Tab_{\\/\mu}$ and $P\in\Tab_{\nu/\emptyset}$
of the standardisations of $T$ and~$\CT{\nu\tr}$ satisfy $K\slid P$ (jeu de
taquin commutes with transposition), so $K\in\Tx$. To show that
$X\in\pi(\F_{\eta,K})$, it suffices to extend the sequence of subspaces
$X_0\ssubset X_1\ssubset\cdots\ssubset X_{\nu_0}$ by interpolation to some
$f\in\F_{\eta,K}$. Now $\eta$ acts as~$0$ on each of the quotient spaces
$X_{i+1}/X_i$, so any choice of complete flags in those spaces leads to a
$f\in\F^{(m)}_\eta$, which has moreover the property that all
partitions~$\mu^i$ occur in the chain~$r_\eta(f)$; we only need to show that
it is possible to obtain $r_\eta(f)=K$. In fact, among the irreducible set of
choices for~$f$, a dense subset has $r_\eta(f)=K$; this follows from the
observation that for any subspace $V'\supset\im\eta$, the projective space
$S=\setof H\in\P^*(V)_\eta: H\supset V'\endset$ meets $U^*_j(\eta)$ whenever
the vertical strip $Y(\\)\setminus Y(J(\eta|_{V'}))$ meets column~$j$, and
then of course the intersection is dense in~$S$ for the minimal such~$j$.

It remains to show that $T$ is a Littlewood-Richardson tableau.
\iftrue
That~$T$ is a semistandard tableau means that each $Y(\mu^{i+1})\setminus
Y(\mu^i)$ is a vertical strip, or equivalently
$(\mu^{i+1})\tr_c\geq(\mu^i)\tr_c\geq(\mu^{i+1})\tr_{c+1}$ for $i,c\in\N$;
this follows from the easily verified inclusions $W_c(\eta|_{X_{i+1}})\supset
W_c(\eta|_{X_i})\supset W_{c+1}(\eta|_{X_{i+1}})$.
\else
One has $X_{i+1}\supset X_i\supset\im(\eta|_{X_{i+1}})$, so it follows from
proposition~\Jsubsetprop\ that $Y(\mu^{i+1})\setminus Y(\mu^i)$ is contained
in the vertical strip $Y(\mu^{i+1})\setminus Y((\mu^{i+1})^\leftarrow)$, where
$(\mu^{i+1})^\leftarrow$ is the Jordan type of the restriction of~$\eta$
to~$\im(\eta|_{X_{i+1}})$ (which is obtained by subtracting~$1$ from every
non-zero part of~$\mu^{i+1}$); therefore $Y(\mu^{i+1})\setminus Y(\mu^i)$ is
itself a vertical strip, and $T$ is semistandard.
\fi
The remaining conditions for~$T$ to be a Littlewood-Richardson tableau can be
formulated in several equivalent ways, but the following will be practical
here, in view of proposition~\Jordanprop: denoting by $T(i,c)$ the number of
entries~$i$ in the first~$c$ rows of~$T$, one has $T(i+1,c+1)\leq T(i,c)$ for
$i,c\in\N$ (this implies $T(i,c)=0$ when $i\geq c$). Now $T(i,c)$ is the
difference between the number of squares in the first~$c$ columns of
$Y(\mu^{i+1})$ and of~$Y(\mu^i)$; putting
$X_i^c=\ker(\eta|_{X_i})^c=X_i\thru\ker\eta^c$, we therefore have by
proposition~\Jordanprop\ that $T(i,c)=\dim X_{i+1}^c-\dim
X_i^c=\dim(X_{i+1}^c/X_i^c)$. Now for all $i,c$ one has
$\eta\inv(X_i^c)=X_{i+1}^{c+1}$, whence $\eta$ induces an injective map
$X_{i+2}^{c+1}/X_{i+1}^{c+1}\to X_{i+1}^c/X_i^c$, giving the required
inequality $T(i+1,c+1)\leq T(i,c)$.
\QED

 \par
\newsection Relative positions of flags.

Besides the use mentioned above of the Robinson-Schensted algorithm as a
computational aid in dealing with classes of jeu de taquin equivalence and
dual equivalence, there is also a direct geometric interpretation, due to
Steinberg, of the correspondence defined by it. To formulate it, we need to
attach a geometric meaning to permutations; it will be based on the fact that
permutations of~$n$ parametrise the orbits for the diagonal action of~$\GL_n$
on~$\F\times\F$. This parametrisation can be defined by associating to each
pair~$(f,f')$ of flags a permutation called the \newterm{relative position} of
$f$~and~$f'$ that will characterise the $\GL_n$-orbit of~$(f,f')$. By
definition $\sigma\in\Sn$ is the relative position of $f$~and~$f'$ if there
exists a basis $\li(e0..n-1)$ of~$V$ such that $f_i=\<\li(e0..i-1)>$ and
$f'_i=\<\li(e\sigma_0..\sigma_{i-1})>$ for all $i\leq n$. The fact that there
always exists a unique such $\sigma$ is the essence of Bruhat's lemma
for~$\GL_n$, but it is useful to give here an explicit demonstration of this
fact. We shall use the auxiliary concept of a growth matrix.

\proclaim \De..
A growth matrix of order~$n$ is a matrix $A=(A_{i,j})_{0\leq i,j\leq n}$ with
entries in~$\N$ satisfying
\item{(i)} $A_{i,0}=A_{0,i}=0$ and $A_{i,n}=A_{n,i}=i$ for $0\leq i\leq n$,
\item{(ii)} $A_{i+1,j}-A_{i,j}\in\{0,1\}$ and $A_{i,j+1}-A_{i,j}\in\{0,1\}$
     for $0\leq i,j<n$,
\item{(iii)} $A_{i+1,j+1}=A_{i,j+1}\implies A_{i+1,j}=A_{i,j}$
     (equivalently $A_{i+1,j+1}=A_{i+1,j}\implies A_{i,j+1}=A_{i,j}$)
     for $0\leq i,j<n$.

A growth matrix $A=(A_{i,j})_{0\leq i,j\leq n}$ corresponds bijectively to a
permutation $\sigma(A)\in\Sn$, whose permutation matrix $\Pi$ (with
$\Pi_{i,j}=\delta_{i,\sigma(A)_j}$ for $0\leq i,j<n$) is related to~$A$ by the
equivalent relations
$$ 
\vcenter{\openup\jot\ialign{\strut\hfil$#$&&$\displaystyle{}#$\hfil\cr
\Pi_{i,j}&=A_{i+1,j+1}-A_{i+1,j}-A_{i,j+1}+A_{i,j}
         &\ttex{for $0\leq i,j<n$,}\cr
A_{i,j}&=\sum_{0\leq i'<i}\ \sum_{0\leq j'<j} \Pi_{i',j'}
       &\ttex{for $0\leq i,j\leq n$.}\cr
}}
$$

Define $\pi(f,f')=\sigma(A)$ where $A$ is the growth matrix with
$A_{i,j}=\dim(f_i\thru f'_j)$. From the definition it follows that if $\sigma$
is the relative position of $f$~and~$f'$, then $\pi(f,f')=\sigma$; in
particular $\sigma$ is unique. Conversely, for~$\sigma=\pi(f,f')$, a basis
$\li(e0..n-1)$ witnessing the fact that $\sigma$ is the relative position of
$f$~and~$f'$ can be constructed: for each $i$ put $j=\sigma\inv_i$ (so that
$\Pi_{i,j}=1$) and choose for $e_i$ any vector in the complement of the
subspace $f_i\thru f'_j$ within $f_{i+1}\thru f'_{j+1}$ (by the construction
of~$\sigma$ this subspace is a hyperplane, and equal to both $f_{i+1}\thru
f'_j$ and $f_i\thru f'_{j+1}$). As an example, if $f=f'$, then
$A_{i,j}=\min(i,j)$, whence $\pi(f,f')$ is the identity permutation; at the
other extreme, when $f$~and~$f'$ are in general position one has
$A_{i,j}=\max(0,i+j-n)$, whence $\pi(f,f')=\w$, the order reversing
permutation.

\heading "Remark"
Besides the mentioned growth matrix~$A$, growth matrices $B$, $C$, and~$D$ can
also be associated to~$(f,f')$ with $B_{i,j}=\dim(f_i/(f_i\thru f'_{n-j}))$, \
$C_{i,j}=\dim(f'_j/(f_{n-i}\thru f'_j))$, and
$D_{i,j}=\dim(V/(f_{n-i}+f'_{n-j}))$. With $\sigma=\sigma(A)=\pi(f,f')$, it
can easily be verified that $\sigma(B)=\sigma\w$, \ $\sigma(C)=\w\sigma$, and
$\sigma(D)=\w\sigma\w$. Since $\dim(V/(f_{n-i}+f'_{n-j}))=\dim(f^*_i\thru
f'^*_j)$, the last case implies that the relative position the dual flags is
obtained by conjugation by~$\w$: $\pi(f^*,f'^*)=\w\pi(f,f')\w$. It is also
obvious that $\pi(f',f)=\pi(f,f')\inv$.

\heading "Remark"
The Bruhat order~`$\leq$' on~$\Sn$ can be defined in terms the associated
growth matrices: if $\sigma=\sigma(A)$ and $\sigma=\sigma(A')$, then one has
$\sigma\leq\sigma'$ if and only if $A_{i,j}\geq A'_{i,j}$ for all~$i,j$. It
follows that for any $\sigma\in\Sn$, the closure of $\setof (f,f'):
\pi(f,f')=\sigma \endset$ in~$\F\times\F$ is $\setof (f,f'):
\pi(f,f')\leq\sigma \endset$.

 \par
\newsection {Interpretation of the Robinson-Schensted algorithm}.

In this section we shall demonstrate the result, due to Steinberg, that for a
pair of flags generically chosen in irreducible components of~$\Fu$
parametrised by a pair of standard Young tableaux, their relative position is
related to those tableaux by the Robinson-Schensted correspondence. Like for
the interpretation of the \Sch. correspondence, this shall be deduced from a
more detailed statement, that gives an interpretation of all the partitions in
the doubly indexed family describing the algorithm; in the current case that
family is the one of \ref{festschrift,~3.2}. We recall the basic
configurations that can occur.

\proclaim\De.. \RSconfigdef
An arrangement of 4~partitions $\ar(\nu,\mu,\mu',\\)$ is called
\item{(i)}  a configuration of type~$RS1$ if $\nu=\mu=\mu'\prec\\$ and
            $\mu_0=\\_0-1$,
\item{(ii)} a configuration of type~$RS2$ if $\nu\prec\mu=\mu'\prec\\$, with
            $\mu_{i+1}=\\_{i+1}-1$ and $\nu_i=\mu_i-1$ for some~$i\geq0$,
\item{(iii)}a configuration of type~$RS3$ if $\nu\prec\mu,\mu'\prec\\$ and
            $\mu\neq\mu'$,
\item{(iv)} a configuration of type~$RS0$ if $\nu=\mu\preceq\mu'=\\$ or
            $\nu=\mu'\preceq\mu=\\$.

\proclaim\De..
A family of partitions $(\\^{[i,j]})$, where $i$~and~$j$ each range over an
interval of~$\Z$, is said to be of type~$RS$ if any configuration
$\ar(\\^{[i,j]},\\^{[i,j+1]},\\^{[i+1,j]},\\^{[i+1,j+1]})$
is of one of the types $RS0$--$RS3$.

Families of type~$RS$ were used in the proof of \ref{festschrift,
theorem~3.2.1} to relate a permutation~$\sigma$ and the pair
$(P,Q)=RS(\sigma)$ of tableaux computed from it by the Robinson-Schensted
algorithm. The pair $(P,Q)$ describes the partitions at the boundary of the
family, while $\sigma$ describes the places where a configuration of
type~$RS1$ occurs; specifying either of these suffices to determine the entire
family. More precisely, to a family $(\\^{[i,j]})_{0\leq i,j\leq n}$ of
type~$RS$ with $\\^{[n,n]}=\\$, we associate tableaux $P,Q\in\Tab_\\$ such
that $\ch P=(\\^{[n,n]},\\^{[n-1,n]},\ldots,\\^{[0,n]})$ and $\ch
Q=(\\^{[n,n]},\\^{[n,n-1]},\ldots,\\^{[n,0]})$, and a permutation
$\sigma=\sigma(A)$ where $A$ is the growth matrix with $A_{i,j}=|\\^{[i,j]}|$,
which means that $i=\sigma_j$ if and only if
$\ar(\\^{[i,j]},\\^{[i,j+1]},\\^{[i+1,j]},\\^{[i+1,j+1]})$ is of type~$RS1$.
(In its relation to $P$, $Q$, and~$\sigma$, our family has its indices
interchanged with respect to \ref{festschrift}; it remains of type~$RS$, and
it still establishes the relation~$(P,Q)=RS(\sigma)$.)

\proclaim Theorem. \RSinterp
Let $P,Q\in\Tab_\\$, and $\sigma\in\Sn$ with $(P,Q)=RS(\sigma)$. For all
$(f,f')$ in a dense open subset of $\F_{\eta,P}\times\F_{\eta,Q}$, the family
$J(\eta|_{f_i\thru f'_j})_{0\leq i,j\leq n}$ is of type~$RS$, and in
particular $\pi(f,f')=\sigma$.

\proclaim Lemma. \intersectlemma
Let $\mu,\mu'\prec\\$, and let $\nu$ be determined by the condition that
$\ar(\nu,\mu,\mu',\\)$ is of one of the types $RS1$--$RS3$. For any
$H_0\in\P^*(V)_\eta$ with $J(\eta|_{H_0})=\mu'$, there is a dense open subset
of $\setof H\in\P^*(V)_\eta: J(\eta|_H)=\mu\endset$ on which $J(\eta|_{H\thru
H_0})=\nu$ holds.

\proof
If $\mu\neq\mu'$ the configuration is of type~$RS3$, and it follows from
proposition~\Jsubsetprop\ that one always has $J(\eta|_{H\thru
H_0})=\nu=\mu\thru\mu'$. Assume now that $\mu=\mu'$, and let $\\-\mu=(i,j)$,
so that $H$ ranges over the set $U^*_j(\eta)$, which has dimension~$i$. If
$i=0$ the configuration is of type~$RS1$ and $U^*_j(\eta)=\{H_0\}$, so that
for the unique choice $H=H_0$ one has $J(\eta|_{H\thru H_0})=
J(\eta|_H)=\mu=\nu$. Finally, if $i>0$, the configuration is of type~$RS2$. We
claim that intersection with~$H_0$ defines a surjective morphism
$\overline{U^*_j(\eta)}\setminus\{H_0\}\to\overline{U^*_{j'}(\eta|_{H_0})}$,
where $\mu-\nu=(i-1,j')$: the lemma then follows by taking for the dense open
subset of~$U^*_j(\eta)$ the intersection of~$U^*_j(\eta)$ with the inverse
image of~$U^*_{j'}(\eta|_{H_0})$. To prove the claim, observe that restriction
to~$H_0$ defines a linear map $W^*_j(\eta)\to
W^*_j(\eta|_{H_0})=W^*_{j'}(\eta|_{H_0})$, whose surjectivity follows from
lemma~\imkerlemma, or from a dimension consideration.
\QED

We note that in the final case each of the fibres of the map
$\overline{U^*_j(\eta)}\setminus\{H_0\}\to\overline{U^*_{j'}(\eta|_{H_0})}$
meets~$U^*_j(\eta)$, whence the restriction of that map to
$U^*_j(\eta)\setminus\{H_0\}$ is still surjective. This means that in this
case any partition obtained from~$\mu$ by removing a corner in some
column~$c>j'$ can arise as $J(\eta|_{H\thru H_0})$ for some non-generic $H\in
U^*_j(\eta)$; of course, by taking $H=H_0$ one can obtain $J(\eta|_{H\thru
H_0})=\mu$ as well.

The lemma can also be formulated as follows: with $\\-\mu=(i,j)$ and
$\mu'-\nu=(i',j')$, one has for $H$ in a dense subset of $U^*_j(\eta)$ that
$H\thru H_0\in U^*_{j'}(\eta|_{H_0})$. It is not however true in general that
as $H$ traverses this subset of $U^*_j(\eta)$, the values $H\thru H_0$
traverse a dense subset of $U^*_{j'}(\eta|_{H_0})$. The reason for this is
that within $H_0$ equipped with $\eta|_{H_0}$ the subspace $\im\eta$ is in no
way special (it is not fixed under automorphisms), yet $H\thru H_0$ always
contains it (because both $H_0$~and~$H$ do so). Although in some cases it can
be shown that all hyperplanes in $U^*_{j'}(\eta|_{H_0})$ contain~$\im\eta$,
there are also cases where this is not so. This circumstance makes the proof
of theorem~\RSinterp\ more difficult than that of theorem~\Schinterp. There
are however no difficulties in extending the lemma as follows.

\proclaim Lemma. \flagcutlemma
Fix an arbitrary flag~$f'\in\Fu$ and $P\in\Tab_\\$; for all~$f$ in a dense
open subset of $\F_{\eta,P}$, the subfamily of
$J(\eta|_{f_i\thru f'_j})_{0\leq i,j\leq n}$ determined by $j\in\{n-1,n\}$ is
of type~$RS$.

Clearly only the hyperplane part $H=f'_{n-1}$ of~$f'$ is relevant here. The
sequence of subspaces $f_i\thru H$ forms a complete $\eta|_H$-stable flag
in~$H$, except that one part is repeated. Denoting this flag (without the
repetition) by $f\thru H$, the lemma states that for generic $f\in\F_{\eta,P}$
one has $f\thru H\in\F_{\eta|_H,P'}$ and $f_i\thru H=f_{i+1}\thru H$, where
$P'$ and~$i$ are found by applying Schensted extraction (reverse insertion)
to~$P$, starting at the square $\\-J(\eta|_H)$ (recall that the entries of~$P$
start at~$0$).

\proof
We apply induction on~$n$. Put $(\\^n,\ldots,\\^0)=\ch P$ (so that
$\\^i=J(\eta|_{f_i})$), and $\mu^i=J(\eta|_{f_i\thru H})$. Applying
lemma~\intersectlemma\ with $\mu'=\mu^n$, $\mu=\\^{n-1}$, and $H_0=H$, we get
for $f_{n-1}$ in a dense open subset~$\U$ of~$\omega(\F_{\eta,P})$ that
$\ar(\mu^{n-1},\\^{n-1},\mu^n,\\^n)$ is of one of the types $RS1$--$RS3$. If
it is of type $RS1$, then $f_{n-1}=H$ for all $f\in\F_{\eta,P}$, so that
$\mu^i=\\^i$ for $i<n$, and the remaining configurations
$\ar(\mu^i,\\^i,\mu^{i+1},\\^{i+1})$ ($0\leq i<n-1$) are of type~$RS0$.
Otherwise we fix any~$K\in\U$ and consider the fibre $\omega\inv(K)$; the
remaining configurations are then taken care of by induction applied to
$\eta|_K$ and $K\thru H$ in place of $\eta$ and $H$, respectively.
\QED

This lemma does not suffice as an induction step for the proof of
theorem~\RSinterp. The induction hypothesis will describe the generic relative
position of the pair formed by a flag in $\F_{\eta|_H,P'}$ and $f'^-$, but for
reasons indicated above, the set of flags $f\thru H$ is not generally dense in
$\F_{\eta|_H,P'}$, and it is conceivable that the generic relative position of
$f\thru H$ and $f'^-$ is different (smaller in the Bruhat order). In order to
dispell this possibility, we shall show that it is possible to write generic
flags $f\in\F_{\eta|_H,P'}$ as $\hat f\thru H$ for some $\hat
f\in\F_{\hat\eta,P}$, where $\hat\eta$ is another nilpotent transformation
of~$V$ than~$\eta$, but with $J(\hat\eta)=\\$ and $\hat\eta|_H=\eta|_H$. To
that end, we shall now consider the reverse process of restricting a nilpotent
transformation to a hyperplane. Let $0\leq i_0\leq n$ and define the following
vector bundle over~$\F_{\eta,T}$:
$$
  \F_{\eta,T,i_0}=\setof (f,v)\in\F_{\eta,T}\times V: v\in f_{i_0}\endset.
$$
Clearly $\F_{\eta,T,i_0}$ is stable under the diagonal $Z_u$-action on
$\F_{\eta,T}\times V$. For any $(f,v)\in\F_{\eta,T,i_0}$ we define in the
$n+1$-dimensional space $V\times k$ a nilpotent transformation~$\hat\eta$ and
a flag~$\hat f$, as follows. Let $e=(0_V,1)\in V\times k$, then $\hat\eta$ is
determined by $\hat\eta|_V=\eta$ and $\hat\eta(e)=v$, and $\hat f$ is defined
by
$$ \hat f_i=\cases{ f_i               & if $0\leq i\leq i_0$, \cr
                    f_{i-1}\oplus\<e> & if $i_0\leq i-1\leq n$. \cr}
$$
It follows from $v\in f_{i_0}$ that $\hat f$ is $\hat\eta$-stable.

\proclaim Lemma. \liftlemma
Let $T\in\Tab_\\$ and $0\leq i_0\leq n$. For all $(f,v)$ in a dense open
subset of $\F_{\eta,T,i_0}$, the associated $\hat\eta$~and~$\hat f$ satisfy
the following property: with $\nu^i=J(\hat\eta|_{\hat f_i})$ and
$\\^i=J(\eta|_{\hat f_i\thru V})$ for $0\leq i\leq n+1$, any configuration of
the form $\ar(\\^i,\nu^i,\\^{i+1},\nu^{i+1})$ is of type~$RS0$ if~$i<i_0$, of
type~$RS1$ if~$i=i_0$, and of type $RS2$~or~$RS3$ if~$i>i_0$.

\proof
One has $\ch T=(\\^{n+1},\ldots,\\^{i_0+1}=\\^{i_0},\ldots,\\^0)$, since $\hat
f_i\thru V$ is $f_i$ if $i\leq i_0$, and $f_{i-1}$ otherwise. Then the unique
sequence of partitions~$\nu^i$ for which the given conditions on
configurations hold is readily constructed; we shall take this as definition
of~$\nu^i$, and show that $J(\hat\eta|_{\hat f_i})=\nu^i$ holds for $0\leq
i\leq n+1$, provided $(f,v)$ is chosen in an appropriate set. The statement
for~$i<i_0$ is immediate. For the remaining statements we shall apply
induction on~$n$ (fixing~$i_0$), starting at $n=i_0$. Put $\nu=\nu^{n+1}$ and
let $\nu-\\=(r,c)$; then by lemma~\imetalemma, the statement $J(\hat\eta)=\nu$
means that $c$ is the minimal value with $v\in\im\eta+\ker\eta^c$. For $n=i_0$
one has $\F_{\eta,T,i_0}=\F_{\eta,T}\times V$, \ $r=0$, and $c=\nu\tr_0$,
which is the least value for which $\im\eta+\ker\eta^c=V$ (in fact
$\ker\eta^c=V$), so the subset of $\F_{\eta,T,i_0}$ of pairs~$(f,v)$ for which
$J(\hat\eta)=\nu$ is dense and open.

As preparation for the induction step, note that from the statement of the
lemma implies that the smallest subspace of~$V$, necessarily $Z_u$-stable,
containing $\im\eta$ and all values $v$ for generic (and hence for all)
$(f,v)\in\F_{\eta,T,i_0}$, equals $\im\eta+\ker\eta^c$. Now assume $n>i_0$,
and choose an arbitrary hyperplane $H\in\omega(\F_{\eta,T})$. Defining
$\hat\omega\:\F_{\eta,T,i_0}\to\P^*(V)$ by $\hat\omega(f,v)=\omega(f)$, it
will suffice to show that the property in the lemma holds in a dense open
subset of $\hat\omega\inv(H)\isom\F_{\eta|_H,T^-,i_0}$, as one can then use
the action of~$Z_u$ to extend this subset to one of~$\F_{\eta,T,i_0}$. The
induction hypothesis gives a dense open subset~$\U$ of~$\hat\omega\inv(H)$ on
which $J(\hat\eta|_{\hat f_i})=\nu^i$ for $0\leq i\leq n$; it remains to show
that generically one also has $J(\hat\eta)=\nu$. If $\nu^n\neq\\$ then this
holds on all of~$\U$ by proposition~\Jsubsetprop. Otherwise
$\ar(\\^n,\nu^n,\\,\nu)$ is of type~$RS2$; let $\nu^n-\\^n=\\-\\^n=(r-1,c')$.
By the induction hypothesis, the smallest subspace of~$H$ containing $\setof
v:(f,v)\in\U\endset$ is $\im(\eta|_H)+\ker(\eta|_H)^{c'}$. Moreover, since
$H\in U^*_{c'}(V)$ this is equal to $\im\eta+\ker\eta^{c'}$ by
lemma~\imkerlemma, and since $c=\\_r$ is the largest part~$\leq c'$ of~$\\$,
also to $\im\eta+\ker\eta^c$. If $c=0$ one has $J(\hat\eta)=\nu$ on all
of~$\U$; otherwise this holds on the open subset
$\setof(f,v)\in\U:v\not\in\im\eta+\ker\eta^{c-1}\endset$, which is non-empty
(since $\\_r=c$ implies
$\dim(\im\eta+\ker\eta^{c-1})<\dim(\im\eta+\ker\eta^c)$) and therefore dense.
\QED

We have now collected all the ingredients necessary to proceed with an
inductive proof of theorem~\RSinterp.

\proofof{theorem~\RSinterp}
For $n<2$ the theorem is obvious, so assume $n\geq2$; we shall prove that for
any choice of $H\in\omega(\F_{\eta,Q})$ the family $J(\eta|_{f_i\thru
f'_j})_{0\leq i,j\leq n}$ is of type~$RS$ for $(f,f')$ in a dense open subset
of $\varphi(H)=\F_{\eta,P}\times(\omega\inv(H)\thru\F_{\eta,Q})$. Let $P'$
and~$i_0$ be the result of applying Schensted extraction to~$P$ starting at
square $\\^{[n,n]}-\\^{[n,n-1]}$, so that $\ch(P')=
(\\^{[n,n-1]},\ldots,\\^{[i_0+1,n-1]}=\\^{[i_0,n-1]},\ldots,\\^{[0,n-1]})$.
Lemma~\flagcutlemma\ provides a dense open subset $\U'$ of~$\varphi(H)$ on
which the configurations
$\ar(\\^{[i-1,n-1]},\\^{[i-1,n]},\\^{[i,n-1]},\\^{[i,n]})$ are all of the
required types, which implies $f\thru H\in\F_{\eta|_H,P'}$. On the other hand
the induction hypothesis provides a dense open subset~$\U''$ of
$\F_{\eta|_H,P'}\times\F_{\eta|_H,Q^-}$ such that remaining configurations
$\ar(\\^{[i-1,j-1]},\\^{[i-1,j]},\\^{[i,j-1]},\\^{[i,j]})$ are of the required
types whenever $(f\thru H,f'^-)\in\U''$. (Here the family of partitions
obtained for $(P',Q^-)$ has been enlarged by duplicating each $\\^{[i_0,j]}$,
to account for the fact that $f_{i_0}\thru H=f_{i_0+1}\thru H$.) The set
$\U=\setof(f,f')\in\U':(f\thru H,f'^-)\in\U''\endset$ is open in~$\varphi(H)$,
and it remains to show that~$\U\neq\emptyset$.

The set~$\U''$ and the dense open set of lemma \liftlemma, applied with
$\eta|_H$ for~$\eta$ and $P'$ for~$T$, both project to a dense open subset of
$\F_{\eta|_H,P'}$, and the intersection of these images is non-empty. For a
choice of inverse images of any point in the intersection, one finds a
nilpotent transformation~$\hat\eta$ of $H\times k$ with $J(\hat\eta)=\\$ and
$\hat\eta|_H=\eta_H$, and $(\hat f,f')\in
\F_{\hat\eta,P}\times\F_{\hat\eta,Q}$ with $J(\hat\eta|_{\hat f_i\thru
f'_j})=\\^{[i,j]}$ for all~$i,j$. There exists a linear \iso.~$g\:H\times k\to
V$ for which $g\after\hat\eta=\eta\after g$, and by the transitive action
of~$Z_u$ on $\omega(\F_{\eta,Q})$ one can achieve that moreover $g(H)=H$; then
$g(\hat f,f')\in\U$, completing the proof.
\QED

\heading "Remark"
Note that in the final argument we do not assume that $g$ fixes~$H$
pointwise; indeed, this may not be possible, since it would imply
$\im\hat\eta=g(\im\hat\eta)=\im\eta$ (the latter identity follows from
$g\after\hat\eta=\eta\after g$), but the whole point of the construction of
$\hat\eta$ was to avoid fixing $\im\hat\eta$ to any particular subspace
of~$H$.

 \par
\vbox{
\newsection {Interpretation of the transposed Robinson-Schensted algorithm}.

In this section we deduce, in close analogy to the previous section, an
interpretation of the version of the Robinson-Schensted algorithm that is defined using
column insertion instead of row insertion, and hence yields transposes of the
tableaux $P$~and~$Q$. We start with a statement that is dual to lemma~\Wlemma.

}

\proclaim Lemma. \Wstarlemma
Let $l\in U_c(\eta)$, and let $\iota\:(V/l)^*\to V^*$ be the map induced by
the canonical projection $V\to V/l$. For $j\neq c$ one has
$\iota(W^*_j(\eta\[l]))=W^*_j(\eta)$, while $\iota(W^*_c(\eta\[l]))$ has
codimension~$1$ in $W^*_c(\eta)$.

\proof
Clearly $\iota$ is injective, and the image of $\im\eta+\ker\eta^j$
under the projection $V\to V/l$ is contained in
$\im\eta\[l]+\ker(\eta\[l])^j$, whence $\iota(W^*_j(\eta\[l]))\subset
W^*_j(\eta)$. The lemma follows by dimension consideration.
\QED

\proclaim\De.. \RStrconfigdef
For any arrangement of partitions $\ar(\nu,\mu,\kappa,\\)$ of type~$RS i$
($i=0,\ldots,3$), the corresponding arrangement of transposed partitions
$\ar(\nu\tr,\mu\tr,\kappa\tr,\\\tr)$ is called a configuration of
type~$RS\tr i$ (the types $RS\tr0$ and $RS\tr3$ are identical to types
$RS0$~and~$RS3$, respectively). If a family of partitions $(\\^{[i,j]})$ is
of type~$RS$, then the family of transposed partitions is said to be of
type~$RS\tr$.

\proclaim Lemma. \modoutlemma
Let $\mu,\mu'\prec\\$, and let $\nu$ be determined by the condition that
$\ar(\nu,\mu,\mu',\\)$ is of one of the types $RS\tr1$--$RS\tr3$. Then for any
$l\in\P(V)_\eta$ with $J(\eta\[l])=\mu'$, there is a dense open subset of
$\setof H\in\P^*(V)_\eta: J(\eta|_H)=\mu\endset$ on which $J(\eta_{H/(l\thru
H)})=\nu$ holds.

\proof
If the configuration is of type~$RS\tr3$, the conclusion follows (without the
need to restrict to a dense subset) by proposition~\Jsubsetprop. Otherwise
$\mu=\mu'$, and we put $\\-\mu=(r,c)$, so that $H$ ranges over the set
$U^*_c(\eta)$ (of dimension~$r$). By lemma~\Wstarlemma, the subspace
$\iota(W^*_c(\eta\[l]))$ of $W^*_c(\eta)$, has codimension~$1$; therefore,
while $H\supset\ker\eta^c$ for all $H\in U^*_c(\eta)$, there is a dense
open subset~$\U$ on which $H\not\supset\eta^{-c}(l)$. If the configuration
is of type~$RS\tr1$ (i.e., $c=0$), one has $l\thru H=\{0\}$ for $H\in\U$,
and therefore $J(\eta_{H/(l\thru H)})=\mu=\nu$. Finally if $c>0$, so that
the configuration is of type~$RS\tr2$, it follows from $H\supset\ker\eta$ or
$l\in\im\eta$ that $H\supset l$; for $H\in\U$ one has
$H/l\not\supset\ker(\eta\[l])^c$, but $H/l\supset\ker(\eta\[l])^{c-1}$ by
application of lemma~\Wstarlemma\ for $j=c-1$, whence
$H/l\in U^*_{c-1}(\eta\[l])$. Then $J(\eta_{H/l})$ differs from
$\mu=J(\eta\[l])$ by a square in column~$c-1$, and is therefore
equal to~$\nu$.
\QED

Note that the generic cases match those of Lemma~\intersectlemma, with all
partitions transposed, but the non-generic cases do not. We know of no
geometric explanation for the first fact, but the latter is no
surprise, since transposition reverses the dominance ordering on Jordan types
that describes orbit closures.

\proclaim Lemma. \flagmodlemma
Let $l\in\P(V)_\eta$. For all flags~$f$ in a dense open subset of $\Fu$, the
following holds: defining $\\^i=J(\eta|_{f_i})$ and
$\mu^i=J(\eta_{f_i/l\thru f_i})$ for $0\leq i\leq n$, any configuration of
the form $\ar(\mu^i,\\^i,\mu^{i+1},\\^{i+1})$ is of one of the
types~$RS\tr0$--$RS\tr3$.

\proof
This follows from lemma~\modoutlemma\ just like lemma~\flagcutlemma\ follows
from \intersectlemma. \QED

Similarly to the situation for the interpretation of the ordinary
Robinson-Schensted correspondence, the sequence of spaces $f_i/(l\thru f_i)$
forms a complete $\eta\[l]$-stable flag in~$V/l$ with one part repeated. For
generic $f\in\F_{\eta,T}$ this flag lies in $\F_{\eta\[l],T'}$ and has
part~$i$ repeated, where this time $T'$~and~$i$ can be found by applying the
Schensted column extraction procedure to~$T$ starting at square
$\\-J(\eta\[l])$, but again the set of flags so obtained is not generally
dense in the indicated set. Here too the difficulty can be resolved by a
converse construction to (in this case) dividing out the line~$l$.

For $0\leq i_0\leq n$ and define the following $Z_u$-stable set (a vector
bundle over~$\F_{\eta,T}$):
$$
  \F_{\eta,T,i_0^*}=
  \setof (f,\phi)\in\F_{\eta,T}\times V^*: f_{i_0}\subset\ker\phi \endset.
$$
To $(f,\phi)\in\F_{\eta,T,i_0^*}$ we associate in the $n+1$-dimensional
space $V\times k$ a nilpotent transformation~$\tilde\eta$ and a flag~$\tilde
f$, as follows. For $(v,x)\in V\times k$ put
$\tilde\eta(v,x)=(\eta(v),\phi(v))$, and with $l=\{0\}\times k\subset V\times
k$ define
$$ \tilde f_i=\cases{ f_i               & if $0\leq i\leq i_0$, \cr
                    f_{i-1}\oplus l   & if $i_0\leq i-1\leq n$. \cr}
$$
Since $f_{i_0}\subset\ker\phi$, the flag $\tilde f$ is
$\tilde\eta$-stable. We identify $(V\times k)/l$ with~$V$, so that
$\tilde f_i/l\isom f_{i-1}$ for $i>i_0$.

\proclaim Lemma. \unmodlemma
Let $T\in\Tab_\\$, and $0\leq i_0\leq n$. For all
$(f,\phi)$ in a $Z_u$-stable dense open subset of $\F_{\eta,T,i_0^*}$, the
associated $\tilde\eta$~and~$\tilde f$ satisfy the following: with
$\nu^i=J(\tilde\eta|_{\tilde f_i})$ and $\\^i=J(\eta_{\tilde f_i/(l\thru
f_i)})$ for $0\leq i\leq n+1$, any configuration of the form
$\ar(\\^i,\nu^i,\\^{i+1},\nu^{i+1})$ is of type~$RS\tr0$ if~$i<i_0$, of
type~$RS\tr1$ if~$i=i_0$, and of type $RS\tr2$~or~$RS3$ if~$i>i_0$.

\proof
The proof is analogous to that of lemma~\liftlemma, with variations as noted
below. Let the families $\nu^i$ and $\\^i$ be determined by the requirements
on the configurations; put $\nu=\nu^{n+1}$ and $\nu-\\=(r,c)$. Since
$(\ker\tilde\eta)/l\isom\ker\eta\thru\ker\phi$, the statement
$J(\tilde\eta)=\nu$ means that $c=\min\setof j:
W_j(\eta)\subset\ker\phi\endset$, by lemma~\keretalemma. In the starting case
$n=i_0$ of the induction this holds, as $c=0$ and $(f,\phi)\in\F_{\eta,T,n^*}$
implies $\phi=0$. For $n>i_0$ and $H\in\omega(\F_{\eta,T})$ there is a
surjection $(f,\phi)\mapsto(f^-,\phi|_H)$ of the fibre $\tilde\omega\inv(H)$
(with $\tilde\omega(f,\phi)=\omega(f)$) onto the variety
$\F_{\eta|_H,T^-,i_0^*}$ (not an isomorphism, as in the proof of
lemma~\liftlemma). The induction hypothesis provides dense open subset of
$\F_{\eta|_H,T^-,i_0^*}$; let $\U\subset\tilde\omega\inv(H)$ be its inverse
image, and $\nu^n-\\^n=(r',c')$, then
$W_{c'}(\eta|_H)=\ker(\eta|_H)\thru\Thru_{(f,\phi)\in\U}\ker(\phi|_H)$. In the
interesting case $RS\tr2$ one has $\nu^n=\\$ and $c'=c+1$, whence
$W_{c'}(\eta)$ strictly contains $W_{c'}(\eta|_H)$ by lemma~\Wlemma, and
$W_{c'}(\eta)\not\subset\ker\phi$ for some $(f,\phi)\in\U$. Therefore
$\setof(f,\phi)\in\U:W_{c'}(\eta)\not\subset\ker\phi\endset$ is dense in~$\U$,
and on this subset $c=\min\setof j:W_j(\eta)\subset\ker\phi\endset$, since
$W_c(\eta)=W_c(\eta|_H)\subset\ker\phi$ for all
$(f,\phi)\in\tilde\omega\inv(H)$.
\QED

\proclaim Theorem. \RStinterp
Let $P,Q\in\Tab_\\$, and $\sigma\in\Sn$ with $(P\tr,Q\tr)=RS(\sigma)$. For all
$(f,f')$ in a dense open subset of $\F_{\eta,P}\times\F^*_{\eta,Q}$, the
family $J(\eta_{f_i/(f_i\thru f'_{n-j})})_{0\leq i,j\leq n}$ is of
type~$RS\tr$, and in particular $\pi(f,f')\w=\sigma$.

\proof
The proof is analogous to that of theorem~\RSinterp, using lemmas
\flagmodlemma~and~\unmodlemma\ instead of \flagcutlemma~and~\liftlemma.
\QED

Note that using theorem~\Schinterp, the equivalence of the descriptions of the
generic value of $\pi(f,f')$ given in theorems \RSinterp~and~\RStinterp\
follows from the purely combinatorial statement \ref{festschrift,
theorem~4.1.1}. We have preferred to give independent proofs of theorems
\RSinterp~and~\RStinterp, thereby giving a geometric explanation for the
``witchcraft operating behind the scenes'' (cf.\ \ref{KnuthACP, p.~60}) of
that combinatorial theorem.

 \par

\Supplement References.

\input \jobname.ref

\bye